\documentclass{amsart}

\usepackage{amssymb}
\usepackage{graphicx}

\newtheorem{theorem}{Theorem}[section]
\newtheorem{lemma}[theorem]{Lemma}
\newtheorem{proposition}[theorem]{Proposition}
\newtheorem{corollary}[theorem]{Corollary}

\newtheorem{_algorithm}[theorem]{Algorithm}

\newtheorem{_definition}[theorem]{Definition}
\newenvironment{definition}{\begin{_definition}\rm}{\end{_definition}}

\newtheorem{_remark}[theorem]{\it Remark}
\newenvironment{remark}{\begin{_remark}\rm}{\end{_remark}}

\newtheorem{_example}[theorem]{Example}

\newtheorem{_assumption}[theorem]{Assumption}

\newtheorem{_construction}[theorem]{Construction}

\newtheorem{_claim}[theorem]{Claim}

\newtheorem{_conjecture}[theorem]{Conjecture}

\newtheorem{_problem}[theorem]{Problem}

\newtheorem{proposition-definition}[theorem]{Proposition-Definition}

\numberwithin{equation}{section}
\numberwithin{table}{section}
\numberwithin{figure}{section}

\newcommand{\C}{\mathord{\mathbb C}}
\newcommand{\F}{\mathord{\mathbb F}}
\renewcommand{\L}{\mathord{\mathbb L}}
\renewcommand{\P}{\mathord{\mathbb  P}}
\newcommand{\Q}{\mathord{\mathbb  Q}}
\newcommand{\R}{\mathord{\mathbb R}}
\newcommand{\Z}{\mathord{\mathbb Z}}

\newcommand{\AAA}{\mathord{\mathcal A}}
\newcommand{\CCC}{\mathord{\mathcal C}}
\newcommand{\DDD}{\mathord{\mathcal D}}
\newcommand{\EEE}{\mathord{\mathcal E}}
\newcommand{\LLL}{\mathord{\mathcal L}}

\newcommand{\OOO}{\mathord{\mathcal O}}
\newcommand{\PPP}{\mathord{\mathcal P}}
\newcommand{\YYY}{\mathord{\mathcal Y}}

\font\mathgot=eufm10

\newcommand{\ppp}{\mathord{\hbox{\mathgot p}}}

\newcommand{\inj}{\hookrightarrow}
\newcommand{\isom}{\mathbin{\,\raise -.6pt\rlap{$\to$}\raise 3.5pt%
\hbox{\hskip .3pt$\mathord{\sim}$}\,}}

\newcommand{\set}[2]{\{\; {#1} \; \mid \; {#2} \;  \}}
\newcommand{\shortset}[2]{\{ {#1} \,|\, {#2}   \}}
\newcommand{\sethd}[3]{\left\{\;\; {#1}\;\;\left|%
\;\;\; \vcenter{\hbox{\parbox{#2}{#3}}}\;\;\;\right.\right\}}

\newcommand{\sprime}{\sp\prime}

\newcommand{\spprime}{\sp{\prime\prime}}
\newcommand{\sptimes}{\sp{\times}}
\newcommand{\sperp}{\sp{\perp}}

\newcommand{\dual}{\sp{\vee}}

\newcommand{\inv}{\sp{-1}}

\newcommand{\Sing}{\operatorname{\rm Sing}\nolimits}
\newcommand{\Hom}{\operatorname{\rm Hom}\nolimits}

\newcommand{\rank}{\operatorname{\rm rank}\nolimits}
\newcommand{\disc}{\operatorname{\rm disc}\nolimits}
\newcommand{\reddisc}{\operatorname{\rm reddisc}\nolimits}
\newcommand{\ord}{\operatorname{\rm ord}\nolimits}
\newcommand{\leng}{\operatorname{\rm leng}\nolimits}
\newcommand{\Pt}{\P^2}

\newcommand{\SL}{\operatorname{\it SL}\nolimits}

\newcommand{\ang}[1]{\langle #1\rangle}

\newcommand{\rmand}{\textrm{and}}
\newcommand{\quand}{\quad\rmand\quad}

\newcommand{\tsum}{\textstyle\sum}

\newcommand{\rootlat}{\Sigma}

\newcommand{\Roots}{\operatorname{\rm Roots}\nolimits}

\newcommand{\erase}[1]{}

\newcommand{\primes}{\PPP}
\newcommand{\pdivs}{\DDD}
\newcommand{\Arith}{\mathord{\rm Arth}}
\newcommand{\Emb}{\mathord{\rm Emb}}
\newcommand{\NK}{\mathord{\rm NK}}
\newcommand{\Klat}{\Lambda_0}
\newcommand{\ssKlat}{\Lambda_{p, \sigma}}
\newcommand{\Zl}{\Z_l}
\newcommand{\Ql}{\Q_l}
\newcommand{\Zt}{\Z_2}

\newcommand{\Zp}{\Z_p}
\newcommand{\Qp}{\Q_p}
\newcommand{\tZ}{\otimes_{\Z}}
\newcommand{\tQ}{\otimes_{\Q}}
\newcommand{\tZl}{\otimes_{\Zl}}

\newcommand{\excess}{\operatorname{\rm-excess}\nolimits}
\newcommand{\lexcess}{l\excess}
\newcommand{\texcess}{2\excess}
\newcommand{\pexcess}{p\excess}
\newcommand{\Sln}{S^{(l)}_n}
\newcommand{\Nln}{N^{(l)}_n}
\newcommand{\spl}{^{(l)}}
\newcommand{\spp}{^{(p)}}
\newcommand{\spt}{^{(2)}}

\newcommand{\LAT}{\L}
\newcommand{\sbps}{_{p, \sigma}}
\newcommand{\sbplus}{_{\mathord{+}}}
\newcommand{\sbminus}{_{\mathord{-}}}
\newcommand{\MW}{\mathord{MW}}
\newcommand{\MWf}{\mathord{MW\hskip -2pt _f}}
\newcommand{\Km}{\mathord{\rm Km}}
\newcommand{\Extlat}{\EEE}
\newcommand{\pol}{H}
\newcommand{\sbXL}{_{(X, \pol)}}
\newcommand{\amp}{\AAA}
\newcommand{\zeroGamma}{{}^0\Gamma}
\newcommand{\Zh}{\Z\hskip .8pt}
\newcommand{\Qh}{\Q\hskip .8pt}
\newcommand{\weight}{\operatorname{\rm wt}\nolimits}

\begin{document}

\title[On normal $K3$ surfaces]{On normal $K3$ surfaces}

\author{Ichiro Shimada}
\address{
Department of Mathematics,
Faculty of Science,
Hokkaido University,
Sapporo 060-0810,
JAPAN
}
\email{shimada@math.sci.hokudai.ac.jp}


\begin{abstract}
We determine all possible configurations  of rational double points
on  complex normal algebraic $K3$ surfaces, 
and  on normal supersingular   $K3$ surfaces
in characteristic $p>19$.
\end{abstract}

\maketitle

\section{Introduction}
In this paper, we mean by  a $K3$ surface an \emph{algebraic} $K3$ surface
defined over an algebraically closed field,
unless otherwise stated.
\par
\medskip
A $K3$ surface $X$ is said to be \emph{supersingular}
(in the sense of Shioda~\cite{MR0374149})
if the rank of the Picard  lattice $S_X$  of $X$ is  $22$.
Supersingular $K3$ surfaces exist only when  the characteristic
of the base field  is positive.
Artin~\cite{MR0371899} showed that,
if $X$ is a supersingular $K3$ surface in characteristic $p>0$,
then the discriminant of  $S_X$  can be 
written as  $-p^{2\sigma_X}$,
where $\sigma_X$ is an integer with $0<\sigma_X\le 10$.
This integer $\sigma_X$ is called the \emph{Artin invariant} of $X$.
\par
\medskip
Let $\Klat$ be an even unimodular $\Zh$-lattice
of rank $22$ with signature $(3, 19)$.
By the structure theorem for unimodular $\Zh$-lattices~(see, for example, Serre~\cite[Chapter V]{MR0344216}),
the $\Zh$-lattice $\Klat$ is unique up to isomorphisms.
If $X$ is a complex $K3$ surface, 
then $H^2 (X, \Z)$ regarded as a $\Z$-lattice by  the cup-product is isomorphic to $\Klat$.
For an \emph{odd} prime integer $p$ and an integer $\sigma$ with $0<\sigma\le 10$,
we denote by 
 $\ssKlat$  an even  $\Zh$-lattice
of rank $22$ with signature $(1, 21)$
such that the discriminant group $\Hom(\ssKlat, \Z)/\ssKlat$
is isomorphic to $(\Z/p\Z)^{\oplus 2\sigma}$.
 Rudakov and Shafarevich~\cite[Theorem in Section 1]{MR633161} showed that
 the $\Zh$-lattice $\ssKlat$ is unique up to isomorphisms.
 If $X$ is a supersingular $K3$ surface 
 in characteristic $p$ with Artin invariant $\sigma$, 
 then  $S_X$ 
 is $p$-elementary by~\cite[Theorem in Section 8]{MR633161}
 and of signature $(1, 21)$ by  the Hodge index theorem,
 and 
 hence 
 $S_X$ is isomorphic to $\ssKlat$.
\par
\medskip
The \emph{primitive closure} of a sublattice $M$ of a $\Zh$-lattice $L$
is  $(M\tZ\Q)\cap L$,
where the intersection is taken in $ L\tZ\Q$.
A sublattice $M\subset  L$ is said to be  \emph{primitive} 
if $(M\tZ\Q)\cap L=M$ holds.
For $\Zh$-lattices $L$ and $L\sprime$, we consider the following condition:
$$
\begin{array}{l@{\;:\;}l}
\Emb (L, L\sprime) &
\parbox[t]{10cm}{There exists a primitive embedding of $L$ into $L\sprime$.} 
\end{array}
$$
We denote by $\primes$ the set of prime integers.
For a non-zero integer $m$,
we denote by $\pdivs (m)\subset \primes$ the set of prime divisors of $m$.
We consider the following arithmetic condition 
on a non-zero integer $d$, a   prime integer $p\in \primes\setminus\pdivs (2d)$, and 
a positive integer $\sigma\le 10$.
$$
\Arith (p, \sigma, d)\;\;\;\;:\;\;\;\; \left(\frac{(-1)^{\sigma+1} d}{p}\right)=-1.
$$
Remark the following:

(i)
Suppose that $d/d\sprime\in (\Q\sptimes)^2$.
Then, for any $p\in \primes\setminus \pdivs(2dd\sprime)$
and any $\sigma$, the conditions
$\Arith (p, \sigma, d)$ and $\Arith (p, \sigma, d\sprime)$ are equivalent.

(ii)
For  fixed $\sigma$ and $d$,
there exists a subset $T_{\sigma, d}$ of $(\Z/4d\Z)\sptimes$ 
such that, 
for $p\in \primes\setminus \pdivs (2d)$,
the condition $\Arith (p, \sigma, d)$ is true
if and only if $p\bmod 4d\in T_{\sigma, d}$.
The set $T_{\sigma, d}$
is empty if and only if
$(-1)^{\sigma+1} d$ is a square integer.
Otherwise, we have $|T_{\sigma, d}|=|(\Z/4d\Z)\sptimes|/2$,
and hence the set of  $p\in \primes\setminus \pdivs (2d)$
for which $\Arith (p, \sigma, d)$ is true
has the natural density $1/2$.
\par
\medskip
The main result of this paper is as follows.
\begin{theorem}\label{thm:emb}
Let $M$ be an even $\Zh$-lattice of rank $r=t\sbplus +t\sbminus $
with signature $(t\sbplus , t\sbminus )$ and of discriminant $d_M$.
Suppose that 
$t\sbplus \le 1$ and $t\sbminus \le 19$.
Then, for a prime integer $p\in \primes\setminus \pdivs (2d_M)$
and a positive integer $\sigma\le 10$, the following hold.
\begin{itemize}
\item[(1)] If $2\sigma> 22-r$, then $\Emb(M, \ssKlat)$ is false.
\item[(2)] If $2\sigma< 22-r$, then $\Emb(M, \ssKlat)$ and $\Emb(M, \Klat)$ are equivalent.
\item[(3)] If $2\sigma= 22-r$, then $\Emb(M, \ssKlat)$ is true if and only if
both of $\Emb(M, \Klat)$ and $\Arith (p, \sigma, d_M)$ are true.
\end{itemize}
\end{theorem}
We  present a geometric  application of Theorem~\ref{thm:emb}.
A \emph{Dynkin type} is a finite formal sum of
 symbols $A_l\,(l\ge 1)$,
$D_m \,(m\ge 4)$ and $E_n \,(n=6,7,8)$
with non-negative  integer coefficients.
For a Dynkin type
$$
R=\sum a_l  A_l +\sum d_m D_m + \sum e_n E_n,
$$
we denote by $\rootlat\sp+_R$ the positive-definite root lattice of type $R$,
and define $\rank (R)$ and $\disc (R)$ to be the rank and the discriminant 
of $\rootlat\sp+_R$:
\begin{eqnarray*}
\rank (R) &:=& \sum a_l  l +\sum d_m m + \sum e_n n,\\
\disc (R) &:=& \prod (l+1)^{a_l}\cdot \prod 4^{d_m} \cdot 3^{e_6}\cdot 2^{e_7}.
\end{eqnarray*}
A \emph{normal $K3$ surface} is a normal surface 
such that its minimal resolution is a $K3$ surface.
It is known that a normal $K3$ surface has only rational double points
as its singularities (Artin~\cite{MR0146182, MR0199191}).
We define the \emph{Dynkin type $R_Y$ of a normal $K3$ surface $Y$}
 to be the Dynkin type of
the singular points on $Y$.
A normal $K3$ surface  is said to be \emph{supersingular}
if its minimal resolution  is supersingular.
The \emph{Artin invariant $\sigma_Y$ of a normal supersingular  $K3$ surface $Y$}
is defined to be the Artin invariant $\sigma_X$ of the minimal resolution $X$ of $Y$.
Note that $\rank (R_Y)$ is equal to the total Milnor number of a normal $K3$ surface $Y$.
In particular,
we have $\rank (R_Y)\le 21$ for any  $Y$,
and $\rank (R_Y)>19$ holds only when $Y$ is supersingular.
\par
\medskip
Let $R$ be a Dynkin type,
$p$ a prime integer, and $\sigma$ a positive integer $\le 10$.
We consider the following  conditions.
$$
\renewcommand{\arraystretch}{1.4}
\begin{array}{l@{\;:\;}l}
\NK (0, R) &
\parbox[t]{10cm}{There exists a complex normal $K3$ surface $Y$
with $R_Y=R$. } \\
\NK (p,\sigma,  R)& 
\parbox[t]{10cm}{There exists a normal supersingular   $K3$ surface $Y$
in characteristic $p$ 
such that $\sigma_Y=\sigma$ and $R_Y=R$.}\\
\NK \sprime(p,\sigma,  R)& 
\parbox[t]{10cm}{Every  supersingular   $K3$ surface $X$
in characteristic $p$ 
with  $\sigma_X=\sigma$ is birational to a normal $K3$ surface $Y$ with  $R_Y=R$.}
\end{array}
$$
We have the following:
\begin{proposition}\label{prop:NN}
The conditions $\NK (p,\sigma,  R)$ and $\NK \sprime(p,\sigma,  R)$ are equivalent.
\end{proposition}
\begin{theorem}\label{thm:NK}
Let $R$ be a Dynkin type with $r:=\rank (R)\le 19$,
and $\sigma$ a positive integer $\le 10$.
We put 
$d_R:=(-1)^r \disc (R)$, and let
 $p$ be an element of $\primes\setminus \pdivs(2d_R)$.
\begin{itemize}
\item[(1)] If $2\sigma> 22-r$, then $\NK(p, \sigma, R)$ is false.
\item[(2)] If $2\sigma< 22-r$, then $\NK(p, \sigma, R)$ and $\NK(0, R)$ are equivalent.
\item[(3)] If $2\sigma= 22-r$, then $\NK(p, \sigma, R)$ is true if and only if
both of $\NK (0, R)$ and $\Arith (p, \sigma, d_R)$ are true.
\end{itemize}
\end{theorem}
For each $p\in \primes$,  a supersingular $K3$ surface 
in characteristic $p$ with Artin invariant $1$
is unique up to isomorphisms~(Ogus~\cite{MR563467, MR717616}).
We denote by  $X_p^{(1)}$  the supersingular $K3$ surface in characteristic $p$
with Artin invariant $1$.
\begin{corollary}
The following  conditions 
on a Dynkin type $R$ with  $r:=\rank (R)\le 19$ are equivalent.
We put $d_R:=(-1)^r \disc (R)$.
\begin{itemize}
\item[(i)] There exists a complex normal $K3$ surface $Y$
with $R_Y=R$.
\item[(ii)] There exists a prime integer  $p\in \primes\setminus\pdivs (2d_R)$ such that  $X_p^{(1)}$
is birational to a normal $K3$ surface $Y$ with $R_Y=R$.
\item[(iii)] For every $p\in \primes\setminus\pdivs (2d_R)$, the supersingular $K3$ surface $X_p^{(1)}$
is birational to a normal $K3$ surface $Y$ with $R_Y=R$.
\end{itemize}
\end{corollary}
Let $Y$ be 
a normal supersingular $K3$ surface in characteristic $p$.
It is proved in~\cite{MR2059747} that, if $\rank (R_Y)=21$, then
$p\in \pdivs(2\disc(R_Y))$ holds.
It is proved in~\cite{Milnor20} that, 
if $\rank (R_Y)=20$, then either $\sigma_Y=1$ or $p\in \pdivs(2\disc(R_Y))$ holds.
(In~\cite{Milnor20}, we have also determined all Dynkin types $R$ of rank $20$
of rational double points that can appear on normal supersingular $K3$ surfaces 
in characteristic $p\notin \pdivs (2 \disc (R))$
with the Artin invariant $1$.)
Therefore, if $\sigma_Y>1$, then either  $\rank (R_Y)\le 19$ or $p\in \pdivs(2\disc(R_Y))$ holds.
Combining this consideration with Theorem~\ref{thm:NK}, 
we obtain restrictions on  Dynkin types of normal supersingular $K3$ surfaces
with large Artin invariants.
\begin{corollary}
Let $Y$ be a normal supersingular $K3$ surface in characteristic $p$
with $\sigma_Y=10$.
Then either one of the following holds.
{\rm (i)} $\rank(R_Y)\le 1$ (that is, $Y$ is smooth or has only one ordinary node as its singularities),  
{\rm (ii)} $R_Y=A_2$ and $p \bmod 24\in \{5, 11, 17, 23\}$, 
{\rm (iii)} $R_Y=2A_1$ and $p \bmod 8\in\{ 3, 7\}$, or {\rm (iv)} $p\in \pdivs(2\disc(R_Y))$.
\end{corollary}
\begin{corollary}
Let $Y$ be a normal supersingular $K3$ surface in characteristic $p$
with $\sigma_Y=9$.
Then either one of the following holds.
{\rm (i)} $\rank(R_Y)\le 3$,
{\rm (ii)} $R_Y=A_4$ and $p\bmod 40\in \{3, 7, 13, 17, 23, 27, 33, 37\}$,  
{\rm (iii)} $R_Y=A_1 + A_3$ and $p \bmod 8\in \{3, 5\}$, 
{\rm (iv)} $R_Y=2 A_1 + A_2$ and $p \bmod 24\in \{5, 7, 17, 19\}$, or {\rm (v)} $p\in \pdivs(2\disc(R_Y))$.
\end{corollary}
Note that, if $p\in \pdivs(2 \disc(R))$ with $\rank(R)\le 21$, then we have $p\le 19$.
Therefore we obtain the following:
\begin{corollary}
The total Milnor number of a normal supersingular $K3$ surface $Y$ in characteristic $p>19$
with Artin invariant $\sigma_Y$ is at most $22-2\sigma_Y$.
\end{corollary}
Let $R$ and $R\sprime$ be Dynkin types.
We write $R\sprime< R$
if the Dynkin diagram of $R\sprime$ can be 
obtained from the Dynkin diagram of $R$
by deleting  some vertexes
and the edges emitting from them.
For a Dynkin type $R$,
we denote by $S(R)$ the set of Dynkin types $R\sprime$ with $R\sprime=R$ or $R\sprime< R$.
A $K3$ surface $X$ is birational to a normal $K3$ surface $Y$
with $R_Y=R$ if and only if there exists a 
configuration of $(-2)$-curves of type $R$ on $X$.
Hence, if $R\sprime\in S(R)$,  we have the following implications:
$$
\NK (0, R)\;\Rightarrow\; \NK (0, R\sprime),
\qquad
\NK (p, \sigma, R)\;\Rightarrow\; \NK (p, \sigma, R\sprime).
$$
We have determined the Boolean value of  $\NK (0, R)$
for each Dynkin type $R$ with $\rank (R) \le 19$, 
and obtained the following:
\begin{theorem}\label{thm:NKlist}
Let $R$ be a Dynkin type of rank $\le 19$.
Then $\NK(0, R)$ is 
true if and only if $S(R)$ does not contain
any  Dynkin type  that appears in Table~\ref{table:List}.
\end{theorem}
\begin{corollary}
Let $R$ be a Dynkin type of rank $\le 14$.
Then there exists a complex normal $K3$ surface $Y$ with
$R_Y=R$.
\end{corollary}
\begin{table}
\parbox{12cm}{
$(\rank  15)$\;
 $A\sb{4} + 11A\sb{1}$,
$2A\sb{2} + 11A\sb{1}$,
$A\sb{2} + 13A\sb{1}$,
\par\medskip
$(\rank  16)$\;
 $3D\sb{4} + 2A\sb{2}$,
$A\sb{6} + A\sb{2} + 8A\sb{1}$,
$A\sb{4} + 2A\sb{2} + 8A\sb{1}$,
\par\medskip
$(\rank  17)$\;
 $E\sb{8} + D\sb{4} + 5A\sb{1}$,
$E\sb{6} + 2D\sb{4} + 3A\sb{1}$,
$E\sb{6} + D\sb{4} + A\sb{2} + 5A\sb{1}$,
$D\sb{7} + 5A\sb{2}$,
$D\sb{5} + 5A\sb{2} + 2A\sb{1}$,
$3D\sb{4} + A\sb{4} + A\sb{1}$,
$2D\sb{4} + A\sb{6} + A\sb{3}$,
$2D\sb{4} + A\sb{6} + 3A\sb{1}$,
$2D\sb{4} + A\sb{4} + A\sb{3} + A\sb{2}$,
$2D\sb{4} + A\sb{4} + A\sb{2} + 3A\sb{1}$,
$2D\sb{4} + 3A\sb{2} + 3A\sb{1}$,
$D\sb{4} + A\sb{8} + 5A\sb{1}$,
$D\sb{4} + 2A\sb{4} + 5A\sb{1}$,
$D\sb{4} + A\sb{3} + 5A\sb{2}$,
$D\sb{4} + 4A\sb{2} + 5A\sb{1}$,
$A\sb{10} + 7A\sb{1}$,
$A\sb{4} + 5A\sb{2} + 3A\sb{1}$,
$A\sb{3} + 5A\sb{2} + 4A\sb{1}$,
$7A\sb{2} + 3A\sb{1}$,
$5A\sb{2} + 7A\sb{1}$,
$17A\sb{1}$,
\par\medskip
$(\rank  18)$\;
 $E\sb{8} + D\sb{4} + 2A\sb{3}$,
$E\sb{6} + D\sb{4} + 2A\sb{3} + A\sb{2}$,
$E\sb{6} + 4A\sb{3}$,
$D\sb{5} + D\sb{4} + 3A\sb{3}$,
$D\sb{4} + A\sb{8} + 2A\sb{3}$,
$D\sb{4} + 2A\sb{4} + 2A\sb{3}$,
$A\sb{7} + 5A\sb{2} + A\sb{1}$,
$2A\sb{4} + 5A\sb{2}$,
$A\sb{4} + 7A\sb{2}$,
$4A\sb{3} + 3A\sb{2}$,
$4A\sb{3} + A\sb{2} + 4A\sb{1}$,
\par\medskip
$(\rank  19)$\;
 $E\sb{7} + 3A\sb{4}$,
$E\sb{7} + 3A\sb{3} + A\sb{2} + A\sb{1}$,
$D\sb{12} + A\sb{7}$,
$D\sb{9} + 3A\sb{3} + A\sb{1}$,
$D\sb{7} + D\sb{5} + 2A\sb{3} + A\sb{1}$,
$D\sb{6} + 2D\sb{5} + A\sb{3}$,
$D\sb{6} + D\sb{5} + 2A\sb{3} + A\sb{2}$,
$D\sb{6} + 3A\sb{4} + A\sb{1}$,
$D\sb{6} + 4A\sb{3} + A\sb{1}$,
$3D\sb{5} + A\sb{3} + A\sb{1}$,
$D\sb{5} + A\sb{5} + 3A\sb{3}$,
$D\sb{5} + 3A\sb{4} + A\sb{2}$,
$D\sb{4} + 4A\sb{3} + 3A\sb{1}$,
$A\sb{7} + 3A\sb{4}$,
$A\sb{6} + 4A\sb{3} + A\sb{1}$,
$A\sb{5} + 3A\sb{4} + A\sb{2}$,
$A\sb{5} + 4A\sb{3} + 2A\sb{1}$,
$A\sb{5} + 3A\sb{3} + 2A\sb{2} + A\sb{1}$,
$3A\sb{4} + 2A\sb{3} + A\sb{1}$,
$3A\sb{4} + A\sb{3} + A\sb{2} + 2A\sb{1}$,
$3A\sb{4} + 2A\sb{2} + 3A\sb{1}$,
$A\sb{4} + 4A\sb{3} + A\sb{2} + A\sb{1}$.
\par\medskip
}
\vskip .3cm
\caption{The  minimal Dynkin types $R$ for which $\NK(0, R)$ is false}
\label{table:List}
\end{table}
Because $p\in \pdivs(2 \disc(R))$ with $\rank(R)\le 21$ implies  $p\le 19$,
Theorems~\ref{thm:NK} and~\ref{thm:NKlist}
combined with the results of
our previous papers~\cite{MR2059747} and~\cite{Milnor20}
determine all  possible configurations  of rational double points
  on normal supersingular   $K3$ surfaces
in characteristic $p>19$.
\par
\medskip
Since $17A_1$ appears in Table~\ref{table:List},
we obtain the following result 
that was proved in  Nikulin~\cite{MR0429917}
for the complex case. See also Section~\ref{subsec:Kummer} of this paper.
\begin{corollary}
{\rm (1)}
There cannot exist 
seventeen disjoint $(-2)$-curves  on a complex $K3$ surface.
{\rm (2)}
There exist seventeen disjoint $(-2)$-curves on a supersingular $K3$ surface
only in characteristic $2$.
\end{corollary}
Remark that,
in characteristic $2$,
there exist twenty-one disjoint $(-2)$-curves on \emph{every}  supersingular $K3$ surface
(\cite{MR2059747, MR2129248}).
\par
\medskip
The proof of Theorems~\ref{thm:emb} and~\ref{thm:NKlist}
is based on the theory of discriminant forms
due to Nikulin~\cite{MR525944},
and the theory of $\lexcess$ due to Conway and Sloane~\cite[Chapter 15]{MR1662447}.
The same method was used in~\cite{MR1813537} to determine
the list  
of   Dynkin types $R_f$ of 
reducible  fibers of 
complex elliptic $K3$ surfaces $f:X\to \P^{1}$ with  a section and the torsion parts $\MWf$ of
their Mordell-Weil groups.
\begin{remark}
Lemma~5.2 in~\cite{MR1813537} is wrong.
It should be replaced with (III) and (IV) in Section~\ref{sec:lexcess}
of the present article.
In the actual calculation of the list of all the pairs $(R_f, \MWf)$
of complex elliptic $K3$ surfaces $f:X\to \P^{1}$ with  a section, 
however,
we  used the correct version of~\cite[Lemma~5.2]{MR1813537}, 
and hence the list  presented in~\cite{MR1813537}
 is valid.
See Remark~\ref{rem:ell}.
\end{remark}
The plan of this paper is as follows.
In Section~\ref{sec:geom},
we prove Proposition~\ref{prop:NN} and deduce Theorem~\ref{thm:NK}
from Theorem~\ref{thm:emb}.
In Section~\ref{sec:lexcess},
we review the theory of $\lexcess$ and 
discriminant forms.
In Section~\ref{sec:proof},
we prove Theorems~\ref{thm:emb} and~\ref{thm:NKlist}.
We conclude the paper with two remarks in the last section.
We give a simple proof of
a theorem of Ogus~\cite[Theorem~7.10]{MR563467} on
supersingular Kummer surfaces,
and investigate, from our point of view,  
the reduction modulo $p$
of a singular $K3$ surface~(in the sense of Shioda and Inose~\cite{MR0441982})
defined over a number field.
\par
\medskip
{\bf Conventions}
\par
(1) Let $D$ be a finite abelian group.
The  \emph{length} of $D$ is the minimal number of generators of $D$,
and is denoted by $\leng (D)$.
\par
(2)
For   $l\in \primes $ and   $x\in \Ql\sptimes$,
we denote by $\ord_l(x)$ the largest integer such that $l^{-\ord_l(x)} x\in \Zl$.
We put $\Z_{\infty}=\Q_{\infty}=\R$.
\par
(3) 
For a divisor $D$ on a $K3$ surface $X$,
let $[D]\in S_X$ denote the class of $D$.
\section{Geometric application}\label{sec:geom}
We prove Proposition~\ref{prop:NN} and deduce Theorem~\ref{thm:NK}
from Theorem~\ref{thm:emb}.
\par
\medskip
Let $X$ be a  $K3$ surface.
A divisor  $\pol$ on  $X$ is called a \emph{polarization}
if $\pol$ is nef, $\pol^2>0$,   and 
the complete linear system  $|\pol|$ has no fixed components.
If $\pol$ is a polarization of  $X$, 
then $|\pol|$ is base-point free by  Saint-Donat~\cite[Corollary 3.2]{MR0364263},
and hence $|\pol|$ defines a morphism $\Phi_{|\pol|}$ from
$X$ to a projective space of dimension $N:=\dim |\pol|=\pol^2/2+1$.
(See Nikulin~\cite[Proposition 0.1]{MR1260944}.)
Let
$$
X\;\;\longrightarrow\;\; Y_{|\pol|}\;\;\longrightarrow\;\;\P^N
$$
be the Stein factorization of $\Phi_{|\pol|}$.
Then $X\to Y_{|\pol|}$ is the  minimal resolution of 
the normal $K3$ surface $Y_{|\pol|}$.
Conversely, let  $X\to Y$ be  the  minimal resolution of 
a normal $K3$ surface $Y$.
Let $\pol\sprime$ be a hyperplane section of $Y$, and  $\pol$  the pull-back of $\pol\sprime$ to $X$.
Then $\pol$ is a polarization of $X$,  and $Y$ is isomorphic to $Y_{|\pol|}$. 
\begin{proposition}\label{prop:polarization}
An element  $v$ of $S_X$ is the  class of a polarization
if and only if $(v,v)>0$, $v$ is nef, and the set
$\shortset{e\in S_X}{(v, e)=1, (e,e)=0}$
is empty.
\end{proposition}
\begin{proof}
See
Nikulin~\cite[Proposition 0.1]{MR1260944},  and 
the argument in the proof of (4)$\Rightarrow$(1)
in Urabe~\cite[Proposition 1.7]{MR1101859}.
\end{proof}
We put
$$
\Xi_X:=\set{v\in S_X}{(v,v)=-2}
\quand
\Gamma_{X} :=\set{x\in S_{X}\tZ\R}{ (x, x)>0}.
$$
For $d\in \Xi_X$, we define the \emph{wall} $d\sperp$ associated with $d$ by 
$$
d\sperp:=\set{x\in S_X\tZ\R}{(x, d)=0}.
$$
Note that the family of walls $d\sperp$ are locally finite in $\Gamma_{X}$.
We denote by 
$$
\zeroGamma_{X} := \set{x\in \Gamma_{X} }{ (x, d)\ne 0\;\;\textrm{for any}\;\; d\in \Xi_{X}}
$$
the complement of these walls in $\Gamma_X$.
Let $W_X$ be 
the subgroup of the orthogonal group $O(S_X)$ of  $S_X$
generated by the reflections
$x \mapsto x+(x, d)d$
into the walls $d\sperp$ associated with the vectors  $d\in \Xi_X$.
Then the subgroup of $O(S_X)$ generated by $W_X$ and $\{\pm 1\}$
acts on the set of connected components of $\zeroGamma_{X}$
transitively.
Let $\amp$ denote the connected component of $\zeroGamma_{X}$
containing the class of a very ample line bundle on $X$.
Then a vector $v\in S_X$ is nef if and only if $v$ is contained in the closure of $\amp$ in $S_{X}\tZ\R$.
Combining these considerations with Proposition~\ref{prop:polarization},
 we obtain the following Corollary. See also~\cite[Proposition 3 in Section 3]{MR633161}.
\begin{corollary}\label{cor:pol}
Let $v\in S_X$ be a vector  such that $(v, v)>0$.
Then there exists an isometry $\phi \in O(S_X)$
such that $\phi(mv)$  is the class of a polarization of $X$
for any integer $m\ge 2$.
\end{corollary}
We introduce a notion from the  lattice theory.
Let $L$ be a negative-definite even $\Zh$-lattice.
A vector $v\in L$ is called a \emph{root} if $(v, v)=-2$ holds.
We denote by $\Roots (L)$ the set of roots in $L$.
A subset $F$ of $\Roots(L)$
is called a \emph{fundamental system of roots in $L$}
if $F$ is a basis of the  sublattice $\ang{\Roots(L)}\subset L$ generated by $\Roots(L)$ and 
each  root $v\in \Roots(L)$ is  written as a linear combination 
$v=\sum_{d\in F}  k_d d$ 
of elements $d$ of $F$ with integer coefficients $k_d$ all non-positive or all non-negative.
Let $t: L\to \R$ be a linear form such that $t(d)\ne 0$
for any $d\in \Roots(L)$.
We put
$$
(\Roots(L))_t^+:=\set{d\in \Roots(L)}{t(d)>0}.
$$
An element $d\in (\Roots(L))_t^+$
is said to be \emph{decomposable}
if there exist vectors $d_1, d_2\in (\Roots(L))_t^+$ such that $d=d_1+d_2$;
otherwise,  we call $d$  \emph{indecomposable}.
The following proposition is proved, for example, in Ebeling~\cite[Proposition 1.4]{MR1938666}.
\begin{proposition-definition}\label{prop:t_roots}
The set $F_t$ of indecomposable elements in $(\Roots(L))_t^+$ is 
a fundamental system of roots in $L$.
We  call  $F_t$
the  \emph{fundamental system of roots associated with $t: L \to \R$}.
\end{proposition-definition}
Let    $\pol$ be  a polarization  of a $K3$ surface $X$.
 The orthogonal complement 
$\ang{[\pol]}\sperp$
of $\ang{[\pol]}$ in $S_X$ is a  negative-definite even lattice.
We put
$$
\Xi\sbXL:=\Roots(\ang{[\pol]}\sperp)=\ang{[\pol]}\sperp \cap \Xi_X.
$$
We denote by  
$F\sbXL$ 
the set of  classes of $(-2)$-curves 
that are contracted by the birational morphism $X\to Y_{|\pol|}$.
It is obvious that $F\sbXL\subset \Xi\sbXL$.
\begin{proposition}\label{prop:geomroots}
The set $F\sbXL$ is equal to the fundamental system of roots $F_\alpha$
in $\ang{[\pol]}\sperp$
associated with the linear form 
$\ang{[\pol]}\sperp \to \R$
given by $v\mapsto (v, \alpha)$, 
where  $\alpha$ is a vector in the connected component $\amp$ of $\zeroGamma_X$.
\end{proposition}
\begin{proof}
We denote by $(\Xi\sbXL)^+_\alpha$
the set of $d\in \Xi\sbXL$ such that $(d, \alpha)>0$.
By the Riemann-Roch theorem,
an element $d\in \Xi\sbXL$ is contained in $(\Xi\sbXL)^+_\alpha$
if and only if $d$ is effective.
Hence we have $F\sbXL\subset (\Xi\sbXL)^+_\alpha$.
Suppose that $[E]\in F\sbXL$ were decomposable in $(\Xi\sbXL)^+_\alpha$,
where $E$ is a $(-2)$-curve contracted by $X\to Y_{|\pol|}$.
Then there would exist   
$[D_1], [D_2]\in (\Xi\sbXL)^+_\alpha$
with $D_1$ and $D_2$ being effective 
such that $[E]=[D_1]+[D_2]$.
Then we would have $D_1+D_2\in |E|$,
which is absurd.
Therefore $[E]$ is indecomposable in  $(\Xi\sbXL)^+_\alpha$,
and  hence $F\sbXL\subset F_\alpha$ is proved.

Conversely, 
let $[D_1], \dots, [D_m]$ be the elements of $F_\alpha$.
Since $F_\alpha\subset (\Xi\sbXL)^+_\alpha$,
we can assume that  $D_1, \dots, D_m$ are effective.
We will show that each $D_i$ is a $(-2)$-curve contracted by $X\to Y_{|\pol|}$.
Let $D_i=F_i+M_i$ be the decomposition of $D_i$
into  the sum of the fixed part $F_i$ and the movable part $M_i$.
Since $\pol$ is nef and $D_i \pol=0$,
we have
$F_i\pol= 0$ and $M_i\pol= 0$.
In particular,
$[M_i]$ is contained in the negative-definite $\Zh$-lattice
$\ang{[\pol]}\sperp$.
Therefore $M_i\ne 0$ would imply  $M_i^2<0$,
which contradicts the movability of $M_i$.
Hence we have  $D_i=F_i$.
Consequently,  the  integral  components $E_1, \dots, E_l$ of  $D_i$
are $(-2)$-curves.
We have
$D_i=a_1 E_1 + \dots + a_l E_l$,
where $a_1, \dots, a_l$ are positive integers.
Since $\pol$ is nef and $D_i\pol=0$,
we have $E_1\pol=\dots =E_l\pol=0$,
and hence   $E_1, \dots, E_l$
are contracted by $\Phi_{|\pol|}$.
Therefore $[E_1], \dots, [E_l]$ are elements
of $F\sbXL\subset F_\alpha$.
Thus, for each $k=1, \dots, l$,
there exists $j_k$ such that $[E_k]=[D_{j_k}]$.
Then we have
$[D_i]=a_1[D_{j_1}]+ \dots +a_l [D_{j_l}]$.
Since 
$[D_1], \dots, [D_m]$
form a basis of the sublattice $\ang{\Xi\sbXL}$ of $\ang{[\pol]}\sperp$,
and $a_1, \dots, a_l$ are positive integers, 
we must have $l=1$, $a_1=1$ and $j_1=i$;
that is, $D_i=E_1$. Hence  $[D_i]\in F\sbXL$ holds, 
and 
$F_\alpha\subset F\sbXL$ is proved.
\end{proof}
\begin{corollary}\label{cor:R}
The Dynkin type of the rational double points on $Y_{|\pol|}$ is 
equal to the Dynkin type of $\Roots(\ang{[\pol]}\sperp)$.
\end{corollary}
Let $L$ be a $\Zh$-lattice.
We denote by $L\dual$ the \emph{dual lattice} $\Hom (L, \Z)$
of $L$.
Then $L$ is embedded in $L\dual$ as a submodule of finite index,
and there exists a natural $\Q$-valued symmetric bilinear form on $L\dual$
that extends  the $\Zh$-valued symmetric bilinear form on $L$.
An \emph{overlattice} of $L$
is a submodule $L\sprime$ of  $L\dual$
containing $L$ such that
the $\Q$-valued  symmetric  bilinear form on $L\dual$
takes values in $\Z$ on $L\sprime$.
If $L$ is embedded in a $\Zh$-lattice $L\spprime$
of the same rank, then $L\spprime$ is naturally embedded in $L\dual$
as an overlattice of $L$.
Let  $L$ be  a negative-definite even $\Zh$-lattice.
If $L\sprime$ is an even overlattice of $L$,
then we have $\Roots (L\sprime)\supseteq\Roots(L)$.
We put 
$$
\Extlat(L):=\sethd{L\sprime}{7cm}{$L\sprime$ is an even overlattice of $L$
 such that $\Roots (L\sprime)=\Roots(L)$ holds}.
$$
For a Dynkin type $R$, we denote by $\rootlat\sp-_R$
the \emph{negative-definite} root lattice of type $R$.
\begin{proposition}\label{prop:RM}
A $K3$ surface $X$ is birational to a normal $K3$ surface $Y$
with $R_Y=R$
if and only if there exists 
 $M\in \Extlat(\rootlat^-_R)$
such that $\Emb (M, S_X)$ is true.
\end{proposition}
\begin{proof}
Combining Corollaries~\ref{cor:pol}~and~\ref{cor:R},
we see that
a $K3$ surface $X$ is birational to a normal $K3$ surface $Y$
with $R_Y=R$
if and only if there exists a vector $v\in S_X$
with $(v, v)>0$ such  that 
$\Roots (\ang{v}\sperp)$  is of type $R$,
where $\ang{v}\sperp$ is the orthogonal complement of  $\ang{v}$ in $S_X$.

Suppose that such a vector $v\in S_X$ exists.
Let $M_0\subset S_X$ be the sublattice of $S_X$ 
generated by $\Roots (\ang{v}\sperp)$.
Then we have an isometry $\varphi: \rootlat^-_R\isom M_0$.
Let $M$ be the overlattice of $\rootlat^-_R$
corresponding by $\varphi$ to the primitive closure of $M_0$ in $S_X$.
Then  $M\in \Extlat(\rootlat^-_R)$ holds, and $\Emb(M, S_X)$ is true.

Conversely,
suppose that there exists  $M\in \Extlat(\rootlat^-_R)$
that admits a primitive embedding $M\inj S_X$.
Let $N$ be the orthogonal complement of $M$ in $S_X$.
Since $M$ is primitive in $S_X$,
the orthogonal complement of $N$ in $S_X$ coincides with $M$.
Hence
a wall $d\sperp$ associated with $d\in \Xi_X$ contains $N\tZ\R$ if and only if 
$d\in \Xi_X\cap M=\Roots (M)=\Roots (\rootlat^-_R)$.
We put
$$
\Gamma_N:=\Gamma_X\cap (N\tZ\R),
$$
which is a non-empty open subset of $N\tZ\R$.
The family of real hyperplanes 
$$
\set{d\sperp\cap (N\tZ\R)}{d\in \Xi_X\setminus \Roots (\rootlat^-_R)}
$$
in $N\tZ\R$ 
is locally finite in $\Gamma_N$, and hence 
there exists   $v\in \Gamma_N\cap N$
such that $v\notin d\sperp$ for any $d\in \Xi_X\setminus \Roots (\rootlat^-_R)$.
Then  $\Roots (\ang{v}\sperp)=\Roots (\rootlat^-_R)$ holds.
\end{proof}
\begin{proposition}\label{prop:N0}
The condition $\NK (0, R)$ is true if and only if
there exists  $M\in \Extlat(\rootlat^-_R)$
such that $\Emb (M, \Klat)$ is true.
\end{proposition}
\begin{proof}
Suppose that there exists a complex normal $K3$ surface $Y$ with $R_Y=R$.
Let $X$ be  the minimal resolution of $Y$.
Then, by Proposition~\ref{prop:RM}, 
there exists $M\in \Extlat(\rootlat^-_R)$
such that $\Emb (M, S_X)$ is true.
Since $S_X$ is primitive in $H^2 (X, \Z)$,
and $H^2 (X, \Z)$ is $\Zh$-isometric  to $\Klat$,
we see that $\Emb (M, \Klat)$ is true.

Conversely, suppose that there exists  $M\in \Extlat(\rootlat^-_R)$
that admits a primitive embedding $M\inj \Klat$.
We choose a vector $h\in \Klat$ such that $(h, h)>0$,
and denote by $S$ the primitive closure of the sublattice of $\Klat$
generated by $M$ and $h$.
Since $M$ is primitive in $\Klat$, the embedding $M\inj S$
is also primitive.
Let $T$ be the orthogonal complement of $S$ in $\Klat$.
We put
$$
\Omega_T:=\set{[\omega]\in \P_* (T\tZ\C)}{(\omega, \omega)=0, (\omega, \bar \omega)>0},
$$
where $[\omega]\subset T\tZ\C$ is the $1$-dimensional linear subspace generated by 
$\omega\in T\tZ\C$.
There exists $[\omega_0]\in \Omega_T$ such that
$\shortset{v\in T}{(\omega_0, v)=0}=\{0\}$.
Then we have
\begin{equation}\label{eq:ST}
\set{v\in \Klat }{(\omega_0, v)=0}=S.
\end{equation}
By the surjectivity of the period mapping for 
complex analytic $K3$ surfaces~(see, for example,~\cite[Chapter VIII]{MR2030225}),
there exist an analytic $K3$ surface $X$ and an isometry
$$
\phi: H^2(X, \Z)\isom \Klat
$$
of $\Zh$-lattices 
such that $\phi\otimes\C$ maps the $1$-dimensional subspace $H^{2,0} (X)\subset H^2(X, \C)$
to  $[\omega_0]$.
By~\eqref{eq:ST}, we have $\phi (S_X)=S$.
Let $h_X\in S_X$ be the vector such that $\phi (h_X)=h$.
Then we have $(h_X, h_X)>0$, and hence $X$ is algebraic.
Since $S$ and $S_X$ is $\Zh$-isometric,
we see that $\Emb (M, S_X)$ is true.
Then $X$ is birational to a normal $K3$ surface $Y$
with $R_Y=R$
by Proposition~\ref{prop:RM}.
\end{proof}
\begin{proof}[Proof of Proposition~\ref{prop:NN} and Theorem~\ref{thm:NK}]
By~\cite[Theorem in Section 8]{MR633161} and~\cite[Theorem in Section 1]{MR633161}
(with~\cite[Proposition in Section 5]{MR633161} for the case of characteristic $2$),
the Picard lattice of a supersingular $K3$ surface is determined,
up to isomorphisms,  by 
the characteristic of the base field and the Artin invariant.
Hence 
Proposition~\ref{prop:NN} follows from   Proposition~\ref{prop:RM}.

Note that $d_R=(-1)^r\disc (R)$ is the discriminant of $\rootlat^-_R$.
If $M$ is an element of $\Extlat(\rootlat^-_R)$
with discriminant $d_M$, then we have $\pdivs (2d_M)\subset \pdivs (2d_R)$, and,
for any  $p\in \primes \setminus \pdivs (2d_R)$,
the conditions
$\Arith (p, \sigma, d_M)$ and $\Arith (p, \sigma, d_R)$
are equivalent,
because  $d_R/d_M =|M/\rootlat^-_R|^2 $ is a square integer.
Therefore Theorem~\ref{thm:NK} follows from Propositions~\ref{prop:RM} and~\ref{prop:N0}  
and Theorem~\ref{thm:emb}.
\end{proof}
\section{The theory of $l$-excess and discriminant forms}\label{sec:lexcess}
See Cassels~\cite{MR522835}, Conway and Sloane~\cite[Chapter 15]{MR1662447}
and Nikulin~\cite{MR525944}
for the details of the results reviewed in this section.
\par
\medskip
Let $R$ be $\Z$, $\Q$, $\Zl$ or  $\Ql$,
where $l\in \primes\cup \{\infty\}$.
An \emph{$R$-lattice} is a free $R$-module $L$ of finite rank
equipped with a non-degenerate symmetric bilinear form
$$
(\phantom{u},\phantom{i})\;:\; L\times L\to R.
$$
We say that $R$-lattices $L$ and $L\sprime$ are \emph{$R$-isometric}
and denote $L\cong L\sprime$
if there exists an isomorphism of $R$-modules $L\isom L\sprime$ that
preserves the symmetric bilinear form.
We sometimes express an $R$-lattice $L$ of rank $n$
by an $n\times n$ symmetric matrix
with components in $R$
by choosing a basis of $L$.
For example, for $a\in R$ with $a\ne 0$,
we denote by $[a]$ the $R$-lattice of rank $1$ generated by
a vector $g$ such that $(g, g)=a$.
For  $R$-lattices $L$ and $L\sprime$,
we denote by $L\oplus L\sprime$ the \emph{orthogonal} direct-sum of $L$ and $L\sprime$.
For $s\in R\setminus\{0\}$,
we denote by $s L$ the $R$-lattice obtained from an $R$-lattice $L$
by multiplying the symmetric bilinear form with $s$.
Suppose that 
an $R$-lattice $L$ is expressed by  a symmetric matrix $M$
with respect to a certain basis of $L$.
Then
$$
\disc (L):=\det(M) \bmod (R\sptimes)^2 \;\;\in\;\; R/(R\sptimes)^2
$$
dose not depend on the choice of the basis of $L$.
We say that $L$ is \emph{unimodular} if $\disc (L)\in R\sptimes/(R\sptimes)^2$.
\par
\medskip
The following is proved in \cite[Theorem 1.2~in Chapter~9]{MR522835}.
\begin{theorem}\label{thm:cassels}
Let $n$ be a positive integer, and $d$ a non-zero integer.
Suppose that, for each $l\in \primes\cup\{\infty\}$,
we are given a $\Zl$-lattice $L_l$
of rank $n$ such that
$\disc (L_l)$ is equal to $d$ 
in $\Zl/(\Zl\sptimes)^2$.
If there exists a $\Qh$-lattice $W$
such that
$W\tQ\Ql$ is $\Ql$-isometric  to $L_l\tZl\Ql$
for each $l\in \primes\cup\{\infty\}$,
then there exists a $\Zh$-lattice $L$ such that
$L\tZ\Zl$ is $\Zl$-isometric  to $L_l$
for each $l\in \primes\cup\{\infty\}$.
\end{theorem}
Let $L$ be an $R$-lattice,
where $R=\Z$ or $\Zl$ with $l\in \primes$,
and let $k$ be the quotient field of $R$.
We put
$$
L\dual :=\Hom _R (L, R).
$$
We have a natural embedding $L\inj L\dual$
of $R$-modules, and a natural $k$-valued symmetric bilinear form
on $L\dual$
that extends the $R$-valued symmetric bilinear form on $L$.
We define the \emph{discriminant group} $D_L$ of $L$ 
by
$$
D_L:=L\dual/L.
$$
If $L$ is a $\Zh$-lattice, then
$\disc (L)=(-1)^{s\sbminus } |D_L|$ holds in $\Z/(\Z\sptimes)^2=\Z$.
\par
\medskip
Suppose that $L$ is a $\Zl$-lattice.
We have an orthogonal direct-sum decomposition
\begin{equation}\label{eq:Jordan}
L=\textstyle\bigoplus_{\nu\ge 0} l^{\nu} L_{\nu},
\end{equation}
where each $L_\nu$ is a unimodular $\Zl$-lattice.
The decomposition~\eqref{eq:Jordan}
is called the \emph{Jordan decomposition} of $L$.
The discriminant group $D_L$ of $L$
is then isomorphic to the direct product 
$\textstyle\prod_{\nu\ge 1} (\Z/l^{\nu}\Z)^{\rank (L_\nu)}$.
In particular, we have
$$
|D_L|=l^{\sum \nu \rank (L_\nu)}
\quand
\leng (D_L)=\rank (L)-\rank (L_0).
$$
We define the \emph{reduced discriminant} of $L$ by 
$$
\reddisc (L):=\textstyle\prod_{\nu\ge 0}\disc (L_\nu)=\disc(L)/|D_L| \;\;\in\;\; \Zl\sptimes/(\Zl\sptimes)^2.
$$

Suppose that $l\ne 2$.
Then
 we have an orthogonal direct-sum decomposition
\begin{equation}\label{eq:oddldecomp}
L\;\;\cong \;\; \textstyle\bigoplus l^{\nu_i} [a_i]\qquad(a_i\in \Z_l\sptimes).
\end{equation}
For $a\in \Zl\sptimes$, we define 
$$
\lexcess (l^{\nu} [a]):=\begin{cases}
(l^\nu -1) \bmod 8 &\textrm{ if $\nu$ is even or $a\in (\Z_l\sptimes)^2$}, \\
(l^\nu +3) \bmod 8 &\textrm{ if $\nu$ is odd and  $a\notin (\Z_l\sptimes)^2$},
\end{cases}
$$
and define $\lexcess (L)\in \Z/8\Z$ to  be the sum of the $l$-excesses of  
the direct summands in~\eqref{eq:oddldecomp}.
It is proved that $\lexcess (L)$ does not  depend on 
the choice of the orthogonal direct-sum decomposition~\eqref{eq:oddldecomp}.
Note that, if $L$ is unimodular, then $\lexcess (L)=0$.

Suppose that $l=2$.
Every unimodular $\Zt$-lattice is $\Zt$-isometric to
an orthogonal direct-sum of copies of the following $\Zt$-lattices:
$$
[a]\quad(a\in \Zt\sptimes), 
\qquad U:=\left[\begin{array}{cc} 0 & 1 \\ 1 & 0 \end{array}\right]\quad\textrm{or}\quad  
 V:=\left[\begin{array}{cc} 2 & 1 \\ 1 & 2 \end{array}\right].
$$
Hence  $L$ has an orthogonal direct-sum decomposition
\begin{equation}\label{eq:twodecomp}
L\;\;\cong\;\;  \textstyle\bigoplus 2\sp{\nu\sb i} [a\sb i] \oplus \textstyle\bigoplus 2\sp{\nu\sb j} U 
\oplus \textstyle\bigoplus 2\sp{\nu\sb k} V, 
\end{equation}
where $a_i\in \Zt\sptimes$.
We put
$$
\renewcommand{\arraystretch}{1.5}
\begin{array}{l}
\texcess (2\sp\nu [a]) :=
\begin{cases}
(1-a) \bmod 8 & \text {\rm if $\nu$ is even or $a\equiv\pm 1 \bmod 8$, } \\
(5-a)  \bmod 8 & \text {\rm if $\nu$ is odd and $a\equiv\pm 3 \bmod 8$, }
\end{cases}
\\
\texcess (2\sp\nu U) := 2  \bmod 8, \qquad
\texcess (2\sp\nu V) := (4 -(-1)^\nu 2) \bmod 8, 
\end{array}
$$
and define $\texcess (L)\in \Z/8\Z$ to  be the sum of the $2$-excesses of  
the direct summands in~\eqref{eq:twodecomp}.
It is proved that $\texcess (L)$ does not  depend on 
the choice of the orthogonal direct-sum decomposition~\eqref{eq:twodecomp}.
The $\texcess$ of a unimodular $\Zt$-lattice need not be $0$.
\par
\medskip
For the following, see Conway and Sloane~\cite[Theorem~8 in Chapter~15]{MR1662447}. 
\begin{theorem}\label{thm:cs}
Let $n$ be a positive integer,
and $d$ a non-zero integer.
Suppose that, for each $l\in \primes\cup\{\infty\}$, we are given a $\Zl$-lattice $L_l$ of rank $n$
such that 
\begin{equation}\label{eq:dc}
\disc (L_l)=d\;\bmod (\Zl\sptimes )^2
\end{equation}
holds in $\Zl/(\Zl\sptimes )^2$. Then there exists a $\Qh$-lattice $W$
such that
$W\tQ\Ql$ is $\Ql$-isometric  to $L_l\tZl\Ql$
for each $l\in \primes\cup\{\infty\}$
if and only if
\begin{equation}\label{eq:pc}
s\sbplus -s\sbminus  +\tsum_{l\in \primes} \lexcess (L_l) \;\equiv\; n\; \bmod 8
\end{equation}
holds, where $(s\sbplus , s\sbminus )$ is the signature of 
the $\R$-lattice $L_{\infty}$.
\end{theorem}
\begin{remark}
If $l\notin \pdivs (2d)$ and $l\ne \infty$, then the condition~\eqref{eq:dc} implies that
the $\Zl$-lattice $L_l$ is  unimodular.
Hence the summation in~\eqref{eq:pc} is in fact finite.
\end{remark}
\begin{definition}
A \emph{finite quadratic form} is a pair $(D, q)$ of a finite abelian group $D$ and 
a map $q: D\to \Q/2\Z$ such that (i)  $q(nx )=n^2 q(x)$ for $n\in \Z$ and $x\in D$, and 
(ii) the map $b: D\times D\to \Q/\Z$ defined by
$b(x, y):=(q(x+y)-q(x)-q(y))/2$ 
is bilinear.
A finite quadratic form $(D, q)$ is said to be \emph{non-degenerate}
if the symmetric bilinear form $b$ is non-degenerate.
\end{definition}
\begin{remark}\label{rem:identify}
Let $(D, q)$ be a finite quadratic form.
Suppose that $D$ is an $l$-group, where $l\in \primes$.
Then the image of $q$ is contained in the subgroup
$$
(\Q/2\Z)_l:=\set{t\in \Q/2\Z}{l^{\nu}t=0\;\; \textrm{for a sufficiently large $\nu$}}=2\Z[1/l]/2\Z
$$
of $\Q/2\Z$. On the other hand,
the canonical homomorphism $\Q/2\Z\to (\Q/2\Z)\tZ\Zl=\Ql/2\Zl$
induces an isomorphism $(\Q/2\Z)_l\isom \Ql/2\Zl$.
Hence we can consider $q$ as a map to $\Ql/2\Zl$.
\erase{
Suppose that $l$ is odd.
Then every element of $\Ql/2\Zl=\Ql/\Zl$
is written uniquely as $a/l^{\nu} \bmod \Zl$,
where $a$ is an integer with $0\le a<l^{\nu}$ and prime to $l$.
The inverse map of the isomorphism $(\Q/2\Z)_l\isom \Ql/2\Zl$ 
is given by
$$
a/l^{\nu} \bmod \Zl\;\;\mapsto\;\;
\begin{cases}
a/l^{\nu} \bmod 2\Z & \textrm{if $a$ is even}, \\
(a+l^{\nu})/l^{\nu} \bmod 2\Z & \textrm{if $a$ is odd}. \\
\end{cases}
$$
}
\end{remark}
\begin{definition}
For a non-degenerate finite quadratic form $(D, q)$ and $l\in \primes$, 
let 
$$
D_l:=\set{t\in D}{l^{\nu}t=0\;\; \textrm{for a sufficiently large $\nu$}}
$$
 denote the $l$-part of $D$, and $q_l$ the restriction
of $q$ to $D_l$. 
We call $(D, q)_l:=(D_l, q_l)$ the \emph{$l$-part} of $(D, q)$.
If $l\notin \pdivs (|D|)$, then  $(D_l, q_l)=(0, 0)$.
We have a decomposition
$$
(D, q)=\textstyle\bigoplus_{l\in \pdivs(|D|)} (D_l, q_l)
$$
that is orthogonal with respect to the symmetric bilinear form $b$.
\end{definition}
Let $R$ be $\Z$ or $\Zl$ with $l\in\primes$, and  $k$ the quotient field  of $R$.
An $R$-lattice $L$ is said to be  \emph{even} if $(v, v)\in 2R$ holds for  every $v\in L$.
Note that, if $l$ is odd, then any $\Zl$-lattice is even.
Note also that a $\Zh$-lattice $L$ is even if and only if the $\Zt$-lattice  $L\tZ\Zt$ is even,
and that 
a $\Zt$-lattice $L$ is even if and only if the  component $L_0$ of 
the Jordan decomposition $L=\bigoplus 2^{\nu}L_{\nu}$ is 
$\Zt$-isometric  to 
an orthogonal direct-sum of copies of $U$ and $V$.
\begin{definition}
For an even $R$-lattice $L$, we can define a map
$$
q_L: D_L\to k/2R
$$
by $q_L(\bar x):=(x, x) \bmod 2R$, where $x\in L\dual$ and $\bar x:=x\bmod L$.
When $R=\Zl$, we consider $q_L$ as a map to $\Q/2\Z$ by 
the isomorphism $\Ql/2\Zl\cong (\Q/2\Z)_l\subset \Q/2\Z$ in Remark~\ref{rem:identify}.
It is easy to see that the finite quadratic form $(D_L, q_L)$
is non-degenerate.
We call $(D_L, q_L)$ the \emph{discriminant form} of $L$.
\end{definition}
We have $\leng(D_L)\le \rank (L)$.
If $L$ is unimodular, then  $(D_L, q_L)=(0,0)$ holds.
If $b_L(\bar x, \bar y):=(q_L(\bar x+\bar y)-q_L(\bar x)-q_L(\bar y))/2$
is the symmetric bilinear form of $(D_L, q_L)$, then
we have $b_L(\bar x, \bar y)=(x, y)\bmod \Z$.
The following is obvious:
\begin{proposition}\label{prop:Datl}
Let $L$ be an even $\Zh$-lattice,
and $l$ a prime integer.
Then the homomorphism
$D_L\to D_{L\tZ\Zl}$
induced from the natural homomorphism
 $L\dual\to L\dual\tZ\Zl=(L\tZ\Zl)\dual$
yields an isomorphism from the $l$-part $(D_L, q_L)_l$
of $(D_L, q_L)$ to 
$(D_{L\tZ\Zl}, q_{L\tZ\Zl})$.
\end{proposition}
Let $(D\spl, q\spl)$ be a non-degenerate  quadratic form
on a finite  abelian $l$-group $D\spl$,
and $n$ a positive integer.
We denote by $\LAT\spl (n, D\spl,  q\spl)$ the set of even $\Zl$-lattices $L$
of rank $n$ such that $(D_L, q_L)$ is isomorphic to $(D\spl, q\spl)$.
We then denote by
$\LLL\spl (n, D\spl, q\spl)\subset\Z/8\Z\times \Zl\sptimes/(\Zl\sptimes)^2$
the image of the map
$$
\renewcommand{\arraystretch}{1.2}
\begin{array}{ccc}
\LAT\spl (n, D\spl,  q\spl)\;\;\; &\to& \Z/8\Z\times \Zl\sptimes/(\Zl\sptimes)^2 \\
L &\mapsto &\;\;\;\tau\spl(L):=[\,\lexcess(L), \reddisc(L)\,].
\end{array}
$$
Let $(D, q)$ be a non-degenerate finite quadratic form, and let 
$$
\LLL\sp{\Z}(n, D, q):=\textstyle\prod_{l\in \pdivs(2|D|)} \LLL\spl (n, D_l, q_l)
$$
be the Cartesian product 
of the sets $\LLL\spl (n, D_l, q_l)$, 
where $(D_l, q_l)$ is the $l$-part of $(D, q)$
and $l$ runs through the prime divisors of $2|D|$.
Let $(s\sbplus , s\sbminus )$ be a pair of non-negative integers such that $s\sbplus +s\sbminus =n$.
We denote by $\LAT\sp{\Z}((s\sbplus , s\sbminus ), D, q)$ the set of even $\Zh$-lattices $L$
of rank $n$ with signature $(s\sbplus , s\sbminus )$ such that 
$(D_L, q_L)$ is isomorphic to $(D, q)$.
By Proposition~\ref{prop:Datl}, we can define a map
$$
\renewcommand{\arraystretch}{1.2}
\begin{array}{ccc}
\LAT\sp{\Z}((s\sbplus , s\sbminus ), D, q)\;\;\; &\to& \LLL\sp{\Z}(n, D, q) \\
L &\mapsto &\;\;\;\tau\sp{\Z}(L):=(\,\tau\spl(L\tZ\Zl) \mid l\in \pdivs(2|D|)\,).
\end{array}
$$
\begin{theorem}\label{thm:mod8}
We put $d:=(-1)^{s\sbminus } |D|$.
Then the image of  $\tau\sp{\Z}$
coincides with the set of
elements
$(\,[\sigma_l, \rho_l]\mid l\in \pdivs(2d)\,)$ of $\LLL\sp{\Z}(n, D, q)$ satisfying the following:
\begin{itemize}
\item[(i)]  $\rho_l= d/ l^{\ord_l(d)}\bmod (\Zl\sptimes)^2$
for each $l\in \pdivs (2d)$, and 
\item[(ii)]
$s\sbplus -s\sbminus  +\sum_{l\in \pdivs (2d)} \sigma_l\;\equiv\; n\mod 8$.
\end{itemize}
In particular, 
the set $\LAT\sp{\Z}((s\sbplus , s\sbminus ), D, q)$ is non-empty if and only if there
exists an element $(\,[\sigma_l, \rho_l]\,|\, l\in \pdivs(2|D|)\,)\in \LLL\sp{\Z}(n, D, q)$
that satisfies  {\rm (i)} and {\rm (ii)}.
\end{theorem}
Let  $l\in \primes$ be an odd prime. 
We choose a non-square  element $v_l\in \Z_l\sptimes$, 
and put $\bar v_l:=v_l\bmod (\Zl\sptimes)^2$,
so that $\Zl\sptimes/(\Zl\sptimes)^2=\{1, \bar v_l\}$ holds. 
We then  define  $\Zl$-lattices $\Sln$ and $\Nln$ of rank $n$ by
\begin{eqnarray*}
\Sln &:=& [1]\oplus\cdots\oplus[1]\oplus[1], \\
\Nln &:=& [1]\oplus\cdots\oplus[1]\oplus[v_l].
\end{eqnarray*}
It is easy to see that  $[v_l]\oplus [v_l]$ is $\Z_l$-isometric  to $[1]\oplus [1]$.
Therefore, 
if $T$  is a unimodular $\Z_l$-lattice of rank $n$, then we have
$$
 T\cong\begin{cases}
 \Sln & \textrm{if $\disc(T)=1$}, \\
 \Nln & \textrm{if $\disc(T)=\bar v_l $}.
 \end{cases}
$$
\begin{proof}[Proof of Theorem~\ref{thm:mod8}]
We denote by $(D_l, q_l)$  the $l$-part of $(D, q)$.
Suppose that $L\in \LAT\sp{\Z}((s\sbplus , s\sbminus ), D, q)$.
Then $\disc (L)=d$ holds.
Since $\disc (L\tZ\Zl)= d \bmod (\Zl\sptimes)^2$
and $|D_{L\tZ\Zl}|=|D_l|=l^{\ord_l (d)}$ by Proposition~\ref{prop:Datl},
we have 
$$
 \reddisc (L\tZ\Zl)\;\;=\;\; d/ l^{\ord_l(d)}\bmod (\Zl\sptimes)^2
$$
for each $l\in \pdivs(2d)$.
Since $\lexcess (L\tZ\Zl)=0$
for every $l\notin \pdivs(2d)$,
we have 
$$
s\sbplus -s\sbminus  +\tsum_{l\in \pdivs (2d)} \lexcess (L\tZ\Zl)\;\;\equiv\;\; n\mod 8
$$
by Theorem~\ref{thm:cs}.
Hence $\tau^{\Z}(L)$ satisfies {\rm (i)} and~{\rm (ii)}.
\par
Conversely, suppose that 
 $([\sigma_l, \rho_l]\,|\, l\in \pdivs(2d))\;\in\; \LLL\sp{\Z}(n, D, q)$
satisfies {\rm (i)} and~{\rm (ii)}.
For each $l\in \pdivs (2d)$, we have an even $\Zl$-lattice $L\spl\in\LAT\spl(n, D_l, q_l)$
such that $\lexcess(L\spl)=\sigma_l$ and $\reddisc (L\spl)=\rho_l$.
Then we have
$$
\disc (L\spl )=\reddisc (L\spl)\cdot |D_l|= d \bmod (\Zl\sptimes)^2
$$
by the condition (i) and $|D_l|=l\sp{\ord_l(d)}$.
For $l\in \primes\setminus \pdivs (2d)$, we put
$$
L\spl:=\begin{cases}
\Sln & \textrm{if $d\in (\Zl\sptimes)^2$},\\
\Nln & \textrm{if $d\notin (\Zl\sptimes)^2$}.
\end{cases}
$$
Then $L\spl\in\LAT\spl(n, D_l, q_l)=\LAT\spl(n, 0, 0)$ 
and $\disc (L\spl )= d \bmod (\Zl\sptimes)^2$ hold.
We put $L^{(\infty)}$ to be an $\R$-lattice of rank $n$
with signature $(s\sbplus , s\sbminus )$.
Then we have 
$\disc (L^{(\infty)})= d \bmod (\R\sptimes)^2$.
Since $\lexcess(L\spl)=0$  for $l\in \primes\setminus \pdivs (2d)$,
the condition (ii) and Theorem~\ref{thm:cs} imply that
there exists a $\Qh$-lattice $W$ of rank $n$
such that $W\tQ\Ql$ is $\Ql$-isometric  to $L\spl\tZl\Ql$
for any $l\in \primes\cup\{\infty\}$.
By Theorem~\ref{thm:cassels},
 there exists a $\Zh$-lattice $L$ of rank $n$
such that $L\tZ\Zl$ is $\Zl$-isometric  to $L\spl$
for any $l\in \primes\cup\{\infty\}$.
Looking at the places $l=2$ and $l=\infty$,
we see that $L$ is  even and of signature $(s\sbplus , s\sbminus )$.
For each $l\in \primes$,
the $l$-part of $(D_L, q_L)$ is 
isomorphic to $(D_{L\spl}, q_{L\spl})\cong (D_l, q_l)$
by Proposition~\ref{prop:Datl}.
Therefore $(D_L, q_L)$ is 
isomorphic to $(D, q)$.
\end{proof}
We fix  $l\in\primes$,
and explain how to calculate the set 
$\LLL\spl (n, D, q)$
for a non-degenerate  quadratic form $(D, q)$
on a finite abelian $l$-group $D$.
\begin{definition}
An orthogonal direct-sum decomposition
$$
(D, q)=(D\sprime, q\sprime)\oplus(D\spprime, q\spprime) 
$$
is said to be \emph{liftable}
if, for any even $\Zl$-lattice $L$
with an isomorphism $\varphi: (D_L, q_L)\isom (D, q)$,
there exists an orthogonal direct-sum decomposition
$L=L\sprime\oplus L\spprime$ such that
$\rank (L\sprime)$ is equal to $\leng (D\sprime)$ and that
$\varphi$ maps $D_{L\sprime}\subset D_L$ to $D\sprime$.
If this is the case,
$\varphi$ induces isomorphisms
$(D_{L\sprime}, q_{L\sprime})\isom (D\sprime, q\sprime)$
and $(D_{L\spprime}, q_{L\spprime})\isom (D\spprime, q\spprime)$.
Therefore we have
$\tau\spl (L\sprime)\in \LLL\spl (\leng(D\sprime), D\sprime, q\sprime)$
and 
$\tau\spl (L\spprime)\in \LLL\spl (n-\leng(D\sprime), D\spprime, q\spprime)$.
\end{definition}
For elements $\tau:=[\sigma, \rho]$ and $\tau\sprime:=[\sigma\sprime,\rho\sprime]$
of $\Z/8\Z\times \Zl\sptimes/(\Zl\sptimes)^2$, we put
$$
\tau*\tau\sprime:=[\sigma+\sigma\sprime, \rho \rho\sprime].
$$
The following is obvious from
$\tau\spl (L\sprime\oplus L\spprime)=\tau\spl (L\sprime)*\tau\spl (L\spprime)$.
\begin{lemma}\label{lem:star}
If an orthogonal direct-sum decomposition $(D, q)=(D\sprime, q\sprime)\oplus(D\spprime, q\spprime)$
is liftable,  then $\LLL\spl (n, D, q)$ is equal to 
$$
\set{\tau*\tau\sprime}{%
\tau \in \LLL\spl (\leng(D\sprime), D\sprime, {q\sprime}),\;
\tau\sprime\in \LLL\spl (n-\leng(D\sprime), D\spprime, q\spprime)}.
$$
\end{lemma}
\begin{lemma}\label{lem:zeroliftable}
The decomposition
$(D, q)=(D, q)\oplus(0, 0)$
is liftable.
\end{lemma}
\begin{proof}
Let  $L$ be an even $\Zl$-lattice
with an isomorphism $(D_L, q_L)\isom (D, q)$,
and let $L=\bigoplus_{\nu\ge 0} l^{\nu} L_{\nu}$ be the Jordan decomposition of $L$.
We put 
$$
L_{\ge 1}:=\textstyle{\bigoplus}_{\nu\ge 1} l^{\nu} L_{\nu}.
$$
Then we have
$\rank  (L_{\ge 1})=\leng(D)$ and $(D_L, q_L)=(D_{L_{\ge 1}}, q_{L_{\ge 1}})$.
Hence the orthogonal direct-sum decomposition $L=L_{\ge 1}\oplus L_0$
has  the required property.
\end{proof}
\begin{lemma}\label{lem:cycliftable}
An orthogonal direct-sum  decomposition
$(D, q)=(D\sprime, q\sprime)\oplus(D\spprime, q\spprime)$
with $D\sprime$ being cyclic 
is liftable.
\end{lemma}
\begin{proof}
Let $l^{\nu}$ be the order of $D\sprime$,
and  $\gamma$  a generator of $D\sprime$.
Since $(D, q)$ is non-degenerate,
so is $(D\sprime, q\sprime)$,
and hence  the order of $b\sprime (\gamma, \gamma)$ in $ \Q/\Z$ is $l^{\nu}$,
where $b\sprime$
is the symmetric bilinear form of $(D\sprime, q\sprime)$.
Let  $L$ be an even $\Zl$-lattice
with an isomorphism $\varphi: (D_L, q_L)\isom (D, q)$.
We choose an element $x\in L\dual$ such that $\varphi(\bar x )=\gamma$,
where $\bar x:= x\bmod L$, 
and put $v:=l^{\nu}x\in L$.
Since $(x, x) \bmod \Zl$ is of order $l^{\nu}$ in $\Ql/\Zl$,
we see that $(v, x)=l^{\nu}(x,x)$ is in $\Zl\sptimes$.
We put $a:=(v, x)\inv \in \Zl\sptimes$.
Since $(w, x)$ is in $\Zl$ and $w- a (w, x) v$ is orthogonal to $v$ for any $w\in L$,
we have  an orthogonal direct-sum decomposition
$L=\ang{v}\oplus\ang{v}\sperp$,
which induces $(D, q)=(D\sprime, q\sprime)\oplus(D\spprime, q\spprime)$
via $\varphi$.
\end{proof}
\begin{definition}
Suppose that $l=2$.
A non-degenerate finite quadratic form $(D, q)$ 
is said to be \emph{of even type}
if $D$ is isomorphic to $\Z/2^{\nu}\Z\times \Z/2^{\nu}\Z$
and the order of $b(\gamma, \gamma)$ in $\Q/\Z$
is strictly smaller than $2^{\nu}$
for any $\gamma\in D$.
\end{definition}
\begin{remark}
Let $L$ be  an even $\Zt$-lattice of rank $2$
with $D_L\cong \Z/2^{\nu}\Z\times \Z/2^{\nu}\Z$.
Then  $(D_L, q_L)$ is of even type
if and only if  $L$ is $\Zt$-isometric to $2^\nu U$ or to $2^\nu V$.
\end{remark}
\begin{lemma}
Suppose that $l=2$.
Then an orthogonal direct-sum  decomposition
$(D, q)=(D\sprime, q\sprime)\oplus(D\spprime, q\spprime)$
with $(D\sprime, q\sprime)$ being of even type  
is liftable.
\end{lemma}
\begin{proof}
Suppose that $D\sprime$ is isomorphic to 
$\Z/2^{\nu}\Z\times \Z/2^{\nu}\Z$,
and let $\gamma_1, \gamma_2$
be elements of $D\sprime$ of order $2^{\nu}$
such that $D\sprime=\ang{\gamma_1}\times \ang{\gamma_2}$.
Since $(D\sprime, q\sprime)$ is of even type,
the orders of $b\sprime (\gamma_1, \gamma_1)$ and $b\sprime (\gamma_2, \gamma_2)$
in $\Q/\Z$ are $< 2^\nu$.
Since $(D\sprime, q\sprime)$ is non-degenerate,
the order of $b\sprime (\gamma_1, \gamma_2)$ 
in $\Q/\Z$ must be   equal to $2^\nu$.
Let  $L$ be an even $\Zt$-lattice
with an isomorphism $\varphi: (D_L, q_L)\isom (D, q)$.
We choose vectors $x_1, x_2\in L\dual$  
such that $\varphi(\bar x_i)=\gamma_i$ for $i=1, 2$,
where $\bar x_i:=x_i\bmod L$,
and put $v_i:=2^{\nu}x_i\in L$.
Then 
there exist $S, T, U\in \Zt$ with $T\in \Zt\sptimes$ such that
$$
\left[
\begin{array}{cc}
(v_1, v_1) & (v_1, v_2) \\
(v_2, v_1) & (v_2, v_2) 
\end{array}
\right]
\;\;=\;\;
2^{\nu}
\left[
\begin{array}{cc}
2 S & T \\ T & 2 U
\end{array}
\right].
$$
Since $4SU-T^2\in \Zt\sptimes$,  the components $\xi_1, \xi_2$ of the vector
$$
\left[
\begin{array}{c}
\xi_1\\
\xi_2
\end{array}
\right]
:=
\left[
\begin{array}{cc}
2 S & T \\ T & 2 U
\end{array}
\right]\inv 
\left[
\begin{array}{c}
(w, x_1)\\
(w, x_2)
\end{array}
\right]
$$
are  elements of $\Zt$
for any $w\in L$.
Moreover,  $w-\xi_1 v_1-\xi_2 v_2$ is orthogonal to 
the sublattice $\ang{v_1, v_2}$ of $L$.
Thus we obtain an orthogonal direct-sum decomposition
$L=\ang{v_1, v_2}\oplus \ang{v_1, v_2}\sperp$,
which induces $(D, q)=(D\sprime, q\sprime)\oplus(D\spprime, q\spprime)$
via $\varphi$.
\end{proof}
\begin{lemma}
If $l$ is odd, then 
$(D, q)$ is an orthogonal direct-sum of finite quadratic forms
on cyclic groups.
If $l=2$, then 
$(D, q)$ is an orthogonal direct-sum of finite quadratic forms
$(D_i, q_i)$,
where, for each $i$, 
$D_i$ is  cyclic or
$(D_i, q_i)$ is  of even type.
\end{lemma}
\begin{proof}
We proceed by induction on $r:=\leng (D)$.
The case where $r=1$ is trivial.
Suppose that $r>1$, and that $D$ is isomorphic to
$\Z/l^{\nu_1}\Z\times \cdots \times \Z/l^{\nu_r}\Z$
with  $\nu_1\ge\dots\ge \nu_r$.
If there exists an element $\gamma\in D$ such that the order of $b(\gamma,\gamma)$ in $\Q/\Z$
is $l^{\nu_1}$,
then $\ang{\gamma}$ is of order $l^{\nu_1}$, and we have an orthogonal direct-sum decomposition
$$
(D, q)=(\ang{\gamma}, q|\ang{\gamma})\oplus (\ang{\gamma}\sperp, q|\ang{\gamma}\sperp)
$$
with $\leng(\ang{\gamma}\sperp)=r-1$.
Suppose that the order of $b(\gamma,\gamma)$ in $\Q/\Z$
is strictly smaller than $l^{\nu_1}$ for any $\gamma\in D$.
Since $(D, q)$ is non-degenerate,
there exist elements $\gamma_1, \gamma_2\in D$ such that
$b(\gamma_1,\gamma_2)\in\Q/\Z$ is of order $l^{\nu_1}$.
If $l\ne 2$, 
then the order of $b(\gamma_1+\gamma_2, \gamma_1+\gamma_2)$ in  $\Q/\Z$ would be  $l^{\nu_1}$.
Therefore we have $l=2$.
We put $D\sprime :=\ang{\gamma_1}\times \ang{\gamma_2}$.
Then $(D\sprime, q|D\sprime)$ is non-degenerate.
We then put
$D\spprime :=D\sp{\prime\perp}$.
Then we have an orthogonal direct-sum decomposition
$$
(D, q)=(D\sprime, q|D\sprime)\oplus (D\spprime , q|D\spprime ),
$$
with $(D\sprime, q|D\sprime)$ being  of even type and $\leng(D\spprime)=r-2$.
\end{proof}
Combining all the results, we can calculate the set $\LLL\spl (n, D, q)$ 
for a positive integer $n$ and a  non-degenerate  quadratic form $(D, q)$
on a finite abelian $l$-group $D$ from the following tables.
\par
\medskip
(I) We have
$$
\LLL\spl (n, D, q)=\emptyset \quad\textrm{if $n< \leng(D)$}.
$$
\par
(II) 
Recall that  $\Zl\sptimes /(\Zl\sptimes)^2=\{1, \bar v_l\}$ for an odd prime $l$.
We also have $\Zt\sptimes /(\Zt\sptimes)^2=\{1,3,5,7\}$.
When $n>0$, we have
\begin{equation*}\label{L0_2}
  \LLL\spl (n, 0, 0)=
\begin{cases}
\{[0, 1], [0, \bar v_l] \} & \textrm{if $l$ is odd,} \\
\emptyset &\textrm{if $l=2$ and $n$ is odd,} \\
\{[n, 1], [n, 5]\} &\textrm{if $l=2$ and $n\equiv 0 \bmod 4$,} \\
\{[n, 3], [n, 7]\}&\textrm{if $l=2$ and $n\equiv 2 \bmod 4$.} \\
\end{cases}
\end{equation*}
\par
(III) {\sl Discriminant forms on cyclic groups.}
Let $\ang{\gamma}$ be a cyclic group of order $l^{\nu}>1$
generated by $\gamma$,
and 
 $q$ a non-degenerate quadratic form on $\ang{\gamma}$.
 Since $q$ is non-degenerate, we can write  $q(\gamma)\in \Q/2\Z$ as $a/l^{\nu}\bmod 2\Z$,
 where $a$ is an integer prime to $l$.
  Suppose that $l$ is odd.
Then we have  
 $$
\LLL\spl (1, \ang{\gamma}, q)=
 \begin{cases}
\{ [ l\sp\nu -1, 1]\} & \text{\rm if $\lambda_l(a)=1$,} \\
\{ [ l\sp\nu -1, \bar v_l]\} & \text{\rm if $\nu$ is even and $\lambda_l(a)=-1$,} \\
\{ [ l\sp\nu +3, \bar v_l]\} & \text{\rm if $\nu$ is odd and $\lambda_l(a)=-1$,} 
\end{cases}
$$
where $\lambda_l: \F_l\sptimes\to\{\pm1\}$ is the Legendre symbol.
When $l=2$,  we have 
$$
\LLL\sp{(2)} (1, \ang{\gamma}, q) =
\begin{cases}
\{ [1-a,a]\} & \text{\rm if $\nu$ is even,} \\
\{ [1-a,a]\} & \text{\rm if $\nu$ is odd, $\nu\ge 2$, and $a\equiv\pm 1\bmod 8$,} \\
\{ [5-a,a]\} & \text{\rm if $\nu$ is odd, $\nu\ge 2$, and $a\equiv\pm 3\bmod 8$,} \\
\{ [0,1],[0,5]\} & \text{\rm if $\nu=1$ and $a\equiv 1\bmod 4$,} \\
\{  [2,3],[2,7]\} & \text{\rm if $\nu=1$ and $a\equiv 3\bmod 4$.} 
\end{cases}
$$
\par
(IV) {\sl Discriminant forms of even type.}
Suppose that $l=2$.
Let  $\ang{\gamma_1}$ and $\ang{\gamma_2}$ be  cyclic groups of order $2^{\nu}$
generated by $\gamma_1$ and $\gamma_2$, where $\nu>0$,
and 
 $q$ a non-degenerate quadratic form on $\ang{\gamma_1}\times \ang{\gamma_2}$
 of even type.
There exist  integers $u$, $v$  and $w$ such that
 $$
  q(\gamma_1)=\frac{2u}{2\sp\nu}\bmod 2\Z,\quad
  q(\gamma_2)=\frac{2w}{2\sp\nu}\bmod 2\Z,\quand
  b(\gamma_1, \gamma_2)=\frac{v}{2\sp\nu}\bmod\Z.
 $$
 Since $q$ is non-degenerate, the integer $v$ is odd.
Then we have
$$
\LLL\sp{(2)} (2, \ang{\gamma_1}\times \ang{\gamma_2}, q) 
=
\begin{cases}
\{ [2, 7]\} & \text{\rm if $uw$ is even,} \\
\{ [2, 3]\} & \text{\rm if $\nu$ is even and $uw$  is odd,} \\
\{ [6, 3]\} & \text{\rm if $\nu$ is odd and $uw$  is odd.} 
\end{cases}
$$
\section{Proof of Main Theorems}\label{sec:proof}
\begin{proposition}\label{prop:genusssKlat}
Let $p$ be an odd prime.
Then $\ssKlat\tZ\Zt$ is $\Zt$-isometric  to $U^{\oplus 11}$,
and $\ssKlat\tZ\Zp$ is $\Zp$-isometric  to
$$
\begin{cases}
S\spp_{22-2\sigma} \oplus p N\spp_{2\sigma} & 
\textrm{if $p\equiv 3 \bmod 4$ and $\sigma \equiv 0 \bmod 2$}, \\
N\spp_{22-2\sigma} \oplus p S\spp_{2\sigma} & 
\textrm{if $p\equiv 3 \bmod 4$ and $\sigma \equiv 1 \bmod 2$}, \\
N\spp_{22-2\sigma} \oplus p N\spp_{2\sigma} & 
\textrm{if $p\equiv 1 \bmod 4$}.
\end{cases}
$$
\end{proposition}
\begin{proof}
Note that $\disc (\ssKlat)=-p^{2\sigma}$.
For simplicity, we put $\Lambda\spl:=\ssKlat\tZ\Zl$.
Since $U\oplus U$ and $V\oplus V$ are $\Zt$-isometric,
the even unimodular $\Zt$-lattice $\Lambda\spt$
is $\Zt$-isometric  to $U^{\oplus 11}$ or to $U^{\oplus 10}\oplus V$.
Since $p^{2\sigma}\in (\Zt\sptimes)^2$, we have $\disc (\Lambda\spt)=-1$
 in $\Zt /(\Zt\sptimes)^2$
 and hence 
 $\Lambda\spt\cong U^{\oplus 11}$.
Therefore we obtain $\texcess(\Lambda\spt)=6$.
Since $D_{\ssKlat}\cong (\Z/p\Z)^{\oplus 2\sigma}$,
the $\Zp$-lattice $\Lambda\spp$ is $\Zp$-isometric  to
$X\oplus pY$, 
where $X$ is either $S\spp_{22-2\sigma} $ or $N\spp_{22-2\sigma} $,
and $Y$ is either $S\spp_{2\sigma} $ or $N\spp_{2\sigma} $.
We have
$$
\pexcess(\Lambda\spp)=
\begin{cases}
2\sigma (p-1)\phantom{+4b} \mod 8 &\textrm{if $Y=S\spp_{2\sigma} $},\\
2\sigma (p-1)+4 \mod 8  &\textrm{if $Y=N\spp_{2\sigma} $}.
\end{cases}
$$
On the other hand, 
from the equality
$$
1-21+\texcess(\Lambda\spt)+\pexcess(\Lambda\spp)\;\;\equiv\;\; 22\mod 8 
$$
in Theorem~\ref{thm:mod8},
we obtain $\pexcess(\Lambda\spp)=4$.
Therefore we have
$$
Y=
\begin{cases}
S\spp_{2\sigma}  &\textrm{if  $2\sigma (p-1)\equiv 4 \bmod 8 $},\\
N\spp_{2\sigma}  &\textrm{if  $2\sigma (p-1)\equiv 0 \bmod 8 $}.\\
\end{cases}
$$
From the equality 
$$
-1=\reddisc (\Lambda\spp)=\disc(X)\disc(Y)=
\begin{cases}
1 &\textrm{if $\disc(X)=\disc(Y)$},\\
\bar v_p  &\textrm{if $\disc(X)\ne \disc(Y)$}
\end{cases}
$$
in $\Zp\sptimes/(\Zp\sptimes)^2$,  we obtain the required result.
\end{proof}
\begin{proposition}\label{prop:Dps}
Let $p$ be an odd prime, and 
let $(D\sbps, q\sbps)$ be the discriminant form of $\ssKlat$.
Then 
$$
\LLL\spp (n, D\sbps, q\sbps)=
\begin{cases}
\emptyset & \textrm{if $n<2\sigma$}, \\
\{[4,1]\} & \textrm{if $n=2\sigma$ and $\sigma (p-1)\equiv 2 \bmod 4$}, \\
\{[4,\bar v_p]\} & \textrm{if $n=2\sigma$ and $\sigma (p-1)\equiv 0 \bmod 4$}, \\
\{[4,1], [4,\bar v_p]\} & \textrm{if $n>2\sigma$}.
\end{cases}
$$
\end{proposition}
\begin{proof}
Let $\ang{\gamma}$ be a cyclic group of order $p$ generated by $\gamma$, and 
let $q_1$ and $q_v$ be the quadratic forms on $\ang{\gamma}$ 
with values in $\Qp/2\Zp=\Qp/\Zp$ 
such that $q_1(\gamma)=1/p\bmod\Zp$ and $q_v(\gamma)=v_p/p\bmod\Zp$,
respectively.
(Let $\tilde v_p\in \Z$ be an integer such that
$\tilde v_p \bmod p=v_p \bmod p\Zp$.
As a quadratic form with values in $\Q/2\Z$,
we have $q_1(\gamma)=(p+1)/p\bmod2\Z$,
 and 
 $$
 q_v(\gamma)=\begin{cases}
 \tilde v_p/p\bmod2\Z & \textrm{  if $\tilde v_p$ is even, } \\
 (\tilde v_p+p)/p\bmod2\Z & \textrm{  if $\tilde v_p$ is odd.}
 \end{cases}
 $$
See Remark~\ref{rem:identify}.)
Then $(\ang{\gamma}, q_1)$ is isomorphic to  the discriminant form  of the $\Zp$-lattice  $p [1]$,
and $(\ang{\gamma}, q_v)$ is isomorphic to  the discriminant form  of the $\Zp$-lattice $p [v_p]$.
By Proposition~\ref{prop:genusssKlat}, we see that 
 $(D\sbps, q\sbps)$ is isomorphic to
 $$
 \begin{cases}
 (\ang{\gamma}, q_1)\sp{\oplus 2\sigma} & \textrm{if $\sigma (p-1)\equiv 2\bmod 4$}, \\
 (\ang{\gamma}, q_1)\sp{\oplus 2\sigma-1} \oplus (\ang{\gamma}, q_v)
 & \textrm{if $\sigma (p-1)\equiv 0\bmod 4$}.
 \end{cases}
 $$
 Hence $\LLL\spp (n, D\sbps, q\sbps)=\emptyset$ for $n<2\sigma$ 
 by (I), and $\LLL\spp (2\sigma, D\sbps, q\sbps)$ is equal to
 $$
  \begin{cases}
 \{[p-1, 1]\sp{*2\sigma} \}=\{[4, 1]\} & \textrm{if $\sigma (p-1)\equiv 2\bmod 4$}, \\
 \{[p-1, 1]\sp{* (2\sigma-1)} * [p+3, \bar v_p] \}=\{[4, \bar v_p]\}
 & \textrm{if $\sigma (p-1)\equiv 0\bmod 4$},
  \end{cases}
 $$
 by Lemmas~\ref{lem:star}~and~\ref{lem:cycliftable} and (III).
If $n>2\sigma$, then $ \LLL\spp (n, D\sbps, q\sbps)$ is equal to 
 $$
\set{\tau*\tau\sprime}{%
\tau\in \LLL\spp (2\sigma, D\sbps, q\sbps), \tau\sprime\in \LLL\spp (n-2\sigma, 0, 0)}
 =\{[4, 1], [4, \bar v_p]\}
 $$
 by Lemmas~\ref{lem:star}~and~\ref{lem:zeroliftable} and (II).
Thus we obtain the required result.
\end{proof}
\begin{proof}[Proof of Theorem~\ref{thm:emb}]
By Nikulin~\cite[Proposition 1.5.1]{MR525944}, 
the condition $\Emb (M, \Klat)$ is true if and only if
\begin{equation}\label{eq:Klat}
\LAT\sp{\Z}((3-t\sbplus, 19-t\sbminus), D_M, -q_M)\ne\emptyset.
\end{equation}
Since $p\notin \pdivs (2d_M)$, 
the condition $\Emb (M, \ssKlat)$ is true if and only if
\begin{equation}\label{eq:ssKlat}
\LAT\sp{\Z}((1-t\sbplus, 21-t\sbminus), D_M\oplus D\sbps, -q_M\oplus q\sbps)\ne\emptyset.
\end{equation}
Remark that 
$$
(-1)^{19-t\sbminus} |D_M|=-d_M
\quand
(-1)^{21-t\sbminus} |D_M\oplus D\sbps|=-p^{2\sigma} d_M.
$$
By Theorem~\ref{thm:mod8}, the condition~\eqref{eq:Klat} is true if and only if there exists
$$
(\,[\sigma_l, \rho_l] \mid l\in \pdivs (2 d_M)\,) \in \LLL\sp{\Z} (22-r, D_M, -q_M)
$$
satisfying
\begin{itemize} 
\setlength{\itemsep}{4pt}
\item[(c1)] $\rho_l= -d_M/l^{\ord_l (d_M)}\bmod (\Zl\sptimes)^2$ for each $l\in \pdivs (2d_M)$, and
\item[(c2)] $-16-t\sbplus+t\sbminus +\sum_{l\in \pdivs(2d_M)} \sigma_l\;\equiv\;22-r \;\mod 8$,
\end{itemize}
and 
the condition~\eqref{eq:ssKlat} is true if and only if there exist
$$
(\,[\sigma\sprime_l, \rho\sprime _l]\,) \in \LLL\sp{\Z} (22-r, D_M, -q_M) 
\quand [\sigma_p, \rho_p]\in \LLL\spp(22-r, D\sbps, q\sbps)
$$
satisfying
\begin{itemize} 
\setlength{\itemsep}{4pt}
\item[(s1)] $\rho\sprime_l= -p^{2\sigma} d_M  /l^{\ord_l(d_M)}\bmod  (\Zl\sptimes)^2$ 
for each $l\in \pdivs (2d_M)$, and 
\item[]  
$\rho_p= -d_M \bmod  (\Zp\sptimes)^2$,  and 
\item[(s2)] $-20-t\sbplus+t\sbminus +\sum_{l\in \pdivs(2d_M)} \sigma\sprime_l+\sigma_p\;\equiv\;22-r \;\mod 8$.
\end{itemize}
Note  that, for $l\in \pdivs (2d_M)$, 
the condition $\rho\sprime_l= -p^{2\sigma}  d_M /l^{\ord_l(d_M)}\bmod  (\Zl\sptimes)^2$
is equivalent to the condition $\rho\sprime_l= -d_M/l^{\ord_l(d_M)}\bmod  (\Zl\sptimes)^2$,
because $p^{2\sigma}\in (\Zl\sptimes)^2$.
By~Proposition~\ref{prop:Dps},
if  $[\sigma_p, \rho_p]\in \LLL\spp(22-r, D\sbps, q\sbps)$, then
$\sigma_p=4$.
Therefore the condition ((s1) and (s2)) is equivalent to the condition 
$$
\parbox{10cm}{(c1) and (c2) and $[4, -d_M]\in \LLL\spp(22-r, D\sbps, q\sbps)$.
}
$$
By~Proposition~\ref{prop:Dps}, 
we have  $[4, -d_M]\in \LLL\spp(22-r, D\sbps, q\sbps)$ 
if and only if $2\sigma<22-r$ holds, or $2\sigma=22-r$ and 
\begin{equation}\label{eq:arthcond}
\begin{array}{l}
(\;\sigma (p-1)\equiv 2 \bmod 4  \quand \lambda _p (-d_M)=1\;)
\;\;\textrm{or}\\    (\;\sigma (p-1)\equiv 0 \bmod 4   \quand  \lambda _p (-d_M)=-1\;)
\end{array}
\end{equation}
hold, 
where $\lambda _p: \F_p\sptimes \to\{\pm 1\}$ is the Legendre symbol.
Since~\eqref{eq:arthcond} is equivalent to $\Arith(p, \sigma, d_M)$, 
Theorem~\ref{thm:emb} is proved.
\end{proof}
\begin{proof}[Proof of Theorem~\ref{thm:NKlist}]
For each Dynkin type $R$ with $r:=\rank (R)\le 19$,
we make the following calculation.
\par
\medskip
(1) We denote by $(D_R, q_R)$ the discriminant form of $\rootlat\sp-_R$,
and by $\Gamma_R$ the image of the natural homomorphism $O(\rootlat\sp-_R)\to O(q_R)$.
(See~\cite[Section 6]{MR1813537} for the description of the group $\Gamma_R$.)
We make the list  of   isotropic subgroups of $(D_R, q_R)$
up to the action of $\Gamma_R$.
By means of Nikulin~\cite[Proposition~1.4.1]{MR525944},
 the list of even overlattices of $\rootlat\sp-_R$ 
 up to the action of $\Gamma_R$ is obtained.
 Then,
 by the method described in~\cite{MR2036331}, 
we make the list $\Extlat(\rootlat\sp-_R)$
up to the action of $\Gamma_R$.

(2) For each $M\in \Extlat(\rootlat\sp-_R)$,
we see whether $\LAT_M:=\LAT\sp{\Z}((3, 19-r), D_M, -q_M)$ is empty or not 
by Theorem~\ref{thm:mod8}.
If we find $M\in \Extlat(\rootlat\sp-_R)$ such that $\LAT_M\ne\emptyset$,
then $\NK (0, R)$ is true.
If  $\LAT_M=\emptyset$ for every $M\in \Extlat(\rootlat\sp-_R)$,
then $\NK (0, R)$ is false.
\end{proof}
\begin{remark}\label{rem:ell}
Let $R$ be a Dynkin type with $r:=\rank (R)\le 18$,
and $\MW$ a finite abelian group.
By~\cite[Theorem 7.1]{MR1813537},
the following are equivalent:
\begin{itemize}
\item[(i)] There exists a complex elliptic $K3$ surface $f: X\to \P\sp 1$
with a section such that the Dynkin type $R_f$ of reducible fibers of $f$ is equal to $R$
and that the torsion part $\MWf$ of the Mordell-Weil  group of $f$ is isomorphic to $\MW$.
\item[(ii)] There exists an element $M\in \Extlat(\rootlat\sp-_R)$
such that $M/\rootlat\sp-_R\cong \MW$ and
that $\LAT\sp{\Z}((2, 18-r), D_M, -q_M)\ne\emptyset$.
\end{itemize}
Therefore, once we have made the list $\Extlat(\rootlat\sp-_R)$
for each Dynkin type $R$ of rank $\le 19$,
it is an easy task to verify the list of all possible  pairs $(R_f, \MWf)$
 given in~\cite{MR1813537}.
\end{remark}
\begin{remark}\label{rem:yang}
Let $\ang{h}$ denote a $\Zh$-lattice of rank $1$ generated by a vector $h$
with $(h, h)=2$.
For  a Dynkin type $R$ with $r:=\rank (R)\le 19$, 
we denote by $\YYY (R)$ the set of even overlattices $M$ of $\rootlat^-_R\oplus \ang{h}$
with the following properties:
\begin{itemize}
\item[(1)] $\Roots (\ang{h}\sperp_M)=\Roots(\rootlat^-_R)$,
where $\ang{h}\sperp_M$ is the orthogonal complement of $\ang{h}$ in $M$, and
\item[(2)] $\shortset{e\in M}{(h, e)=1, (e,e)=0 }=\emptyset$.
\end{itemize}
By Yang~\cite{MR1387816}, 
the following are equivalent:

\begin{itemize}
\item[(i)] There exists a complex reduced  plane curve $C\subset\Pt$ of degree $6$
with only simple singularities such that 
the Dynkin type  of  $\Sing(C)$ is equal to $R$.
\item[(ii)] There exists an element  $M\in \YYY(R)$
such that 
 $\LAT\sp{\Z}((2, 19-r), D_M, -q_M)\ne\emptyset$.
\end{itemize}
During the proof of Theorem~\ref{thm:NKlist},
we have also calculated  the set $\YYY(R)$
for each $R$,
and confirmed the validity of  Yang's list~\cite{MR1387816} of  configurations of singular points of  
 complex sextic curves  with only simple singularities.
\end{remark}
\section{Concluding remarks}\label{sec:rems}
\subsection{Kummer surfaces}\label{subsec:Kummer}
We work over an algebraically closed field of characteristic $p>0$ with $p\ne 2$.
Let $A$ be an abelian surface, and $\iota:A\to A$ the inversion.
Then $Y_A:=A/\ang{\iota}$ is a normal $K3$ surface with $R_{Y_A}=16A_1$.
The minimal resolution $\Km(A)$  of $Y_A$ is called the \emph{Kummer surface}.
We give a simple proof of 
the following theorem due to Ogus~\cite[Theorem~7.10]{MR563467}.
\begin{theorem}\label{thm:Km}
A supersingular $K3$ surface is a Kummer surface
if and only if the Artin invariant is $1$ or $2$.
\end{theorem}
\begin{proof}
Since $\NK (0, 16A_1)$ is true and $\Arith (p, 3, (-1)^{16}2^{16})$ is false,
Theorem~\ref{thm:NK} implies that $\NK (p, \sigma, 16A_1)$ is true
if and only if $\sigma\le 2$.
Thus the ``only if\hskip .5 pt" part of Theorem~\ref{thm:Km} is proved.
To show the ``if\hskip .5 pt" part,
it is enough to prove that the minimal resolution
of a normal $K3$ surface $Y$ with $R_Y=16A_1$
is a Kummer surface.
For this purpose, we use the following Lemma,
which can be  checked easily   by using  a computer:
\begin{lemma}\label{lem:bincode}
Let $\CCC$ be 
a  binary linear code of length $16$ with dimension $\ge 5$
such that the weight $\weight (w)$ of every word $w$ 
satisfies $\weight (w)\equiv 0\bmod 4$ and $\weight (w)\ne 4$.
Then there exists a word of weight $16$ in $\CCC$.
\end{lemma}
We consider subgroups of the discriminant group $D_{16A_1}\cong \F_2\sp{\oplus 16}$ of $\Sigma\sp-_{16A_1}$
as  binary linear codes of length $16$.%
\begin{lemma}\label{lem:lat16A1}
If  $M\in \Extlat (\Sigma\sp-_{16A_1})$ satisfies $\leng(D_M)\le 6$,
then $M/\Sigma\sp-_{16A_1}\subset D_{16A_1}$ contains a word of weight $16$.
\end{lemma}
\begin{proof}
Let $\CCC\subset D_{16A_1}$ be a linear code.
Then $\CCC$ is isotropic with respect to $q_{16A_1}$
if and only if $\weight (w)\equiv 0\bmod 4$  for every $w\in \CCC$.
Suppose that $\CCC$ is isotropic.
Then the corresponding even overlattice $M_{\CCC}$ of $\Sigma\sp-_{16A_1}$ satisfies 
$\Roots (M_{\CCC})=\Roots(\Sigma\sp-_{16A_1})$  if and only if
$\weight (w)\ne 4$ for every $w\in \CCC$.
Because  $\leng (D_{M_{\CCC}})=16-2\dim\CCC$ by  Nikulin~\cite[Proposition~1.4.1]{MR525944},
we obtain Lemma~\ref{lem:lat16A1} from Lemma~\ref{lem:bincode}.
\end{proof}
Suppose that
$Y$ is  a normal $K3$ surface with $R_Y=16A_1$,
and $X\to Y$  the minimal resolution.
We denote by $\Sigma_X$ the sublattice 
of $S_X$ generated by the classes 
of the $(-2)$-curves $E_1, \dots, E_{16}$ contracted by $X\to Y$,
and let  $M_X$ be the primitive closure of $\Sigma_X$ in $S_X$.
Then we have $M_X\in \Extlat (\Sigma_X)$ by Proposition~\ref{prop:geomroots}.
Moreover we have $\leng (D_{M_X})\le 6$,
because $\Emb (M_X, \ssKlat)$ is true,
where $\sigma=\sigma_X$,
and hence $\LLL\spt (22-\rank(M_X), D_{M_X}, -q_{M_X})\ne \emptyset$.
By Lemma~\ref{lem:lat16A1},
there exists a word of weight $16$ in  the code $M_X/\Sigma_X$.
Hence we have $([E_1]+\cdots+[E_{16}])/2\in M_X$.
Therefore there exists a double covering $A\sprime\to X$
whose branch locus is $E_1\cup \dots\cup E_{16}$.
Then the contraction of $(-1)$-curves on $A\sprime$
yields an abelian surface $A$,
and $X$ is isomorphic to the Kummer surface $\Km(A)$. (See~\cite[Lemma 7.12]{MR563467}).
\end{proof}
\begin{remark}
In fact,
a linear  code $\CCC\subset \F_2^{\oplus 16}$  with 
the properties described  in Lemma~\ref{lem:bincode} is unique up to isomorphisms.
See Nikulin~\cite{MR0429917} for the description of this code in terms of 
$4$-dimensional affine geometry over $\F_2$.
\end{remark}
\subsection{Singular $K3$ surfaces}\label{subsec:ShiodaInose}
A complex $K3$ surface $X$ is called \emph{singular} (in the sense of Shioda and Inose~\cite{MR0441982})
if $S_X$ is of rank $20$.
Let $X$ be a singular $K3$ surface,
and  $T_X$ the transcendental lattice of $X$.
Then $T_X$ possesses a canonical orientation $\eta_X$
determined by the holomorphic $2$-form on $X$.
Shioda and Inose~\cite{MR0441982} showed that
the mapping $X\mapsto (T_X, \eta_X)$
induces a bijection from  the set of isomorphism classes
of singular $K3$ surfaces to the set of 
$\SL_2(\Z)$-equivalence classes of positive-definite 
even binary forms.

In~\cite{MR0441982}, it is also shown that 
every  singular $K3$ surface $X$ can be  defined over 
a number field $F$.
(See  Inose~\cite{MR578868} for an explicit defining equation.)
For a maximal ideal $\ppp$ of the integer ring $\OOO_F$ of $F$,
let $X(\ppp)$ denote the reduction of $X$ at $\ppp$.
\begin{proposition}
Suppose that   a singular $K3$ surface $X$ is defined over a number field $F$.
Let $\ppp$ be a maximal ideal  of $\OOO_F$
with  residue  characteristic $p$.
Suppose that 
$p$ is prime to $2\disc (T_X)$,  and that 
$X(\ppp)$ is a supersingular $K3$ surface.
Then the Artin invariant of $X(\ppp)$ is $1$, and we have 
\begin{equation}\label{eq:ssprime}
\left(\frac{-\disc (T_X)}{p}\right)=-1.
\end{equation}
\end{proposition}
\begin{proof}
Since the signature of $S_X$ is $(1, 19)$,
we have $\disc (S_X)=-\disc(T_X)$.
Let $\sigma$ be the Artin invariant of $X(\ppp)$.
The reduction induces an embedding $S_X\inj S_{X(\ppp)}$.
Let $M$ be the primitive closure of $S_X$ in $S_{X(\ppp)}$.
Then $\Emb (M, \ssKlat)$ is true.
Since $M$ is of rank $20$ and $\disc (S_X)/\disc (M)$ is a square integer,
it follows from Theorem~\ref{thm:emb} that  $\sigma=1$,
and  that $\Arith (p, 1, \disc (S_X))$ is true.
Therefore  we obtain~\eqref{eq:ssprime}.
\end{proof}
\bibliographystyle{plain}

\def\cprime{$'$} \def\cprime{$'$} \def\cprime{$'$} \def\cprime{$'$}

\end{document}